\documentclass[a4paper,10pt]{article}

\usepackage[utf8]{inputenc}
\usepackage{amssymb}
\usepackage{amsmath}
\usepackage{graphicx}
\usepackage{xcolor}
\usepackage{natbib}
\usepackage{url}
\usepackage{algorithm}
\usepackage{algpseudocode}
\usepackage{algorithmicx}

\usepackage[left=2cm,top=2.5cm,right=2cm,bottom=2.5cm]{geometry}

\graphicspath{{Figure/}}
\newcommand{\bEp}{\mathbf{E}'}
\newcommand{\bHp}{\mathbf{H}'}
\newcommand{\bJp}{\mathbf{J}'}
\newcommand{\bMp}{\mathbf{M}'}
\newcommand{\Ep}{{E}'}
\newcommand{\Hp}{{H}'}
\newcommand{\Jp}{{J}'}

\newcommand{\bsig}{\mathbf{\sigma}}

\newcommand{\veps}{{\varepsilon}}

\begin{document}

\renewcommand{\thefootnote}{\fnsymbol{footnote}}
\title{Controlled-source electromagnetic modelling using high order finite-difference time-domain method on a nonuniform grid}

\author{Pengliang Yang$^1$ and Rune Mittet$^2$\\
  $^1$School of Mathematics, Harbin Institute of Technology,  Harbin, China, 150001\\
  E-mail: ypl.2100@gmail.com\\
  $^2$Norwegian University of Science and Technology (NTNU), Norway\\
  E-mail: mittet.rune@gmail.com\\
}

\maketitle

\begin{abstract}
  Simulation of 3D low-frequency electromagnetic fields propagating in the Earth is computationally
  expensive.
  We present a fictitious wave domain high-order finite-difference time-domain (FDTD) modelling method on nonuniform grids to compute frequency-domain 3D controlled-source electromagnetic (CSEM) data.
  The method overcomes the inconsistency issue widely present in the conventional 2nd order staggered grid finite difference scheme over nonuniform grid, achieving high accuracy with arbitrarily high order scheme. The finite-difference coefficients adaptive to the node spacings, can be accurately computed by inverting a Vandermonde matrix system using efficient algorithm.
  A generic stability condition applicable to nonuniform grids is established, revealing the dependence of the time step and these finite-difference coefficients.
  A recursion scheme using fixed point iterations is designed to determine the stretching
  factor to generate the optimal nonuniform grid.   The grid stretching in our method reduces the number of grid points required in the discretization, making it more efficient than the
  standard high-order FDTD with a densely sampled uniform grid.
  Instead of stretching in both vertical and horizontal directions,
  better accuracy of our method is observed when the grid is stretched
  along the depth without horizontal stretching.  The efficiency and accuracy of our method are
  demonstrated by numerical examples. 
\end{abstract}

\section{Introduction}

Marine controlled source electromagnetics (CSEM) provides valuable information about
subsurface resistivities and therefore potentially about pore fluids or rocks. It is very useful to decipher
subsurface properties to assist energy exploration, in particular when combined with seismic data.
The CSEM technology relies on low-frequency electromagnetic field propagation to probe the
subsurface. The low-frequency electromagnetic (EM) field propagation does not lend itself to an intuitive
understanding in the same manner as seismic field propagation does due to the diffusive nature
of the EM field in conductive media. Thus, three-dimensional modeling
becomes an important tool for the interpretation of CSEM data. Imaging of marine CSEM data is today
mainly done by inversion of the observed electric and/or magnetic fields.

The kernel of CSEM inversion is the numerical simulation of 3D electromagnetic field propagation, which
is computationally expensive. Reducing the simulation time without compromising the accuracy is important.
It can shorten the turnaround time for an imaging project while reducing the investments in computer hardware.
The implementation of nonuniform grid schemes is a well known strategy to reduce
simulation time. To retain good accuracy we propose a high-order finite-difference approach.

There are many studies on diffusive electromagnetic modelling using
different methods. Examples are the frequency-domain finite-difference method
\citep{newman1995frequency,smith1996conservative1,mulder2006multigrid,streich20093d},
the frequency-domain finite-element method \citep{li20072d,da2012finite,key2016mare2dem,rochlitz2019custem},
and the time-domain finite-difference method \citep{oristaglio1984diffusion,wang1993finite,Taflove_2005_CEF}. 
A key fact in all numerical modelling methods is that the computational cost and
the memory requirement are connected and cannot be splitted. A method can be
very efficient if more computer memory is available. The efficiency and accuracy of the modelling can be dramatically
hampered when the available computer resources are restricted.

Due to the diffusive nature of low-frequency CSEM fields, most of the 3D CSEM modelling schemes resort
to the frequency-domain solution of the Maxwell equation to avoid the high computational cost dictated
by the restrictive stability condition for the direct solution in the time domain. Time-domain methods are
attractive options because they require less amount of memory than frequency-domain modelling
within a model of the same size. Another advantage with
time-domain solutions is that multiple frequencies can be extracted from the same simulation.
Both the frequency-domain finite-difference method and the frequency-domain finite-element method
formulate Maxwell equation as a linear equation system, which may be solved using
direct \citep{streich20093d} or iterative
\citep{smith1996conservative2,mulder2006multigrid,puzyrev2013parallel} solvers.
A nice feature with a direct solver is that multiple right-hand sides are fast to calculate
after the system matrix has been factorized or inverted. However, there are significant implementation
challenges with this approach when realistic size marine CSEM surveys are simulated. The memory requirements
are large even if the equation system is sparse.

The finite-difference time-domain (FDTD) modelling based on the staggered grid
proposed by \citet{Yee_1966_NSI} has for several decades been a main workhorse for many
EM applications. The implementation of the numerical core is straight forward and the 
computational efficiency is good for wave phenomena. The computational efficiency for diffusive
phenomena is rather poor in the time domain. The system of partial differential equations can be considered stiff
in this case \citep{Mittet_2010_HFD} and a very small time step is required to retain stability.
The computational efficiency can be improved significantly due to a
correspondence principle for wave and diffusion fields \citep{lee1989new,de1996general,Mittet_2010_HFD}.
\citet{Maao_2007_FFT} proposed a mixed wave and diffusion-domain
FDTD method to perform numerically efficient CSEM modelling. This method allowed for large time
steps compared to a purely diffusion-domain solution.
\citet{Mittet_2010_HFD} proposed a high-order FDTD scheme by utilizing the fictitious
wave to diffusion-domain transformation. The simulation is performed in the wave domain where
the propagation velocity is proportional to the square root of resistivity.

\begin{figure}
  \centering
  \includegraphics[width=0.75\linewidth]{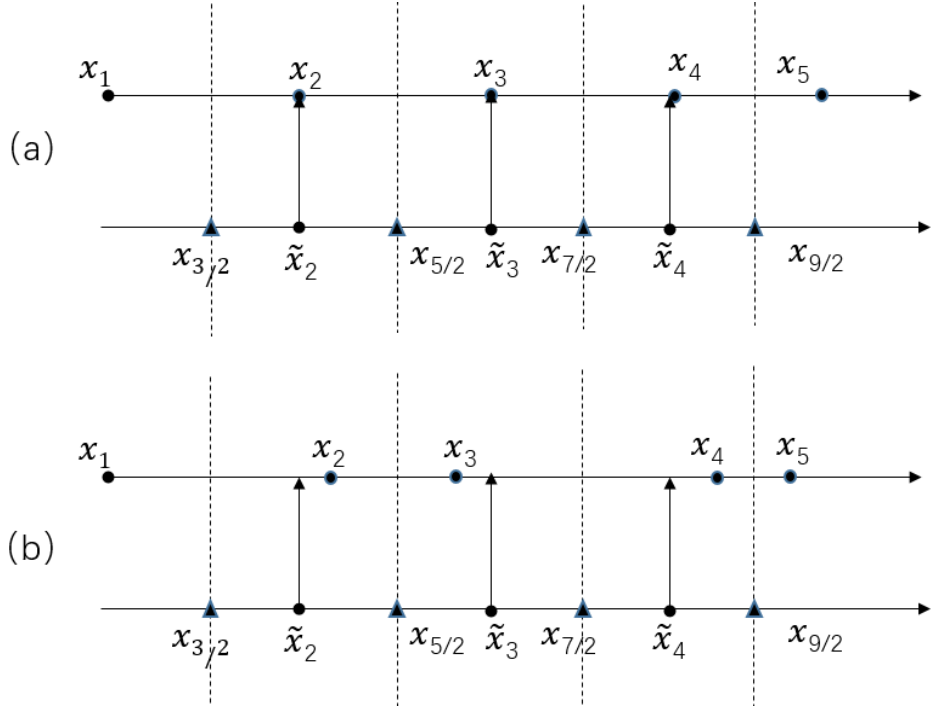}
  \caption{The 2nd order staggered-grid finite-difference scheme on (a) uniform grid and (b) nonuniform
    grid.  The grid points on staggered grid are obtained by taking the midpoints from nonuniform
    grid to ensure the 2nd order accuracy. This accuracy is not guaranteed on nonuniform grid since the
    midpoints is inconsistent. For example, $x_{3/2}= (x_1+x_2)/2$, $x_{5/2}=(x_2+x_3)/2$,
    but $\tilde{x}_2=(x_{3/2} + x_{5/2})/2 \neq x_2$ due to uniform grid spacing.
    This inconsistency leads to first-order local truncation error using staggered-grid
    finite-difference scheme on the nonuniform grid.}\label{fig:figure1}
\end{figure}

The Yee grid (staggered grid) FDTD scheme is often the method of choice due to a good agreement with physics.
It gives divergence free magnetic fields and electric currents \citep{smith1996conservative1}. 
This standard scheme proposed by \citet{Yee_1966_NSI} is based on the second-order approximation
of the first derivatives assuming an equispaced mesh. On the uniform grid, moving from second-order
FDTD to high-order FDTD is straight forward \citep{Mittet_2010_HFD}, and gives improved modelling
accuracy due to the reduction of spatial dispersion errors.
To improve the modelling efficiency, the use of nonuniform grid  is widespread
\citep{newman1995frequency,mulder2006multigrid}, however it results in inconsistencies for the
grid staggering. As illustrated in Figure~\ref{fig:figure1}, this inconsistency leads to
only first-order local truncation error, even though the global accuracy may be up to second order
\citep{monk1994convergence}. The problem is persistent and has remained unresolved for high-order schemes.

We propose an efficient 3D CSEM simulation method with high accuracy using
\emph{high-order  FDTD on a staggered, nonuniform grid}, following the fictitious wave domain approach
\citep{Mittet_2010_HFD}.  To resolve the inconsistency issue in conventional 2nd order staggered grid approach,
our key recognition is that the order of local truncation error can be
arbitrarily high also on a nonuniform grid if the finite-difference operator coefficients are adapted
properly to the variable grid spacing. The derivative operator coefficients are calculated by inverting
a Vandermonde matrix system.

To gain good modeling efficiency, we transforms the diffusive Maxwell equation into
the fictitious wave domain. The efficiency of the method is restricted by stability
condition: the stepsize in time is proportional to the inverse of the propagation velocity of the field,
affected by the node spacing. The gridding of the same physical domain leads to different number of gridpoints,
affecting the size of the linear system to be solved. The use of nonuniform grid helps to
reduce the number of gridpoint, thus reducing the computational cost. Unfortunately,
the stability condition for high order FDTD over non-uniform grid is non-trivial. An important contribution of this paper
is to establish a new stability condition valid for arbitrarily high order FDTD scheme on nonuniform grid.
The stability condition shows the strong dependence between time step and the finite difference
coefficients computed by inverting the Vandermonde matrix. This is not known in EM geophysics community, as far as
we know.

To generate the optimal grid using a power law, we design a recursion scheme
using fixed point iterations to find the optimal stretching factor. We prove the recursion
scheme is guaranteed to converge. The optimal factor
found by the recursive scheme allows accurate matching of the computational domain
using given number of mesh points.
The high accuracy and efficiency of this high-order FDTD method on a nonuniform grid
will be exemplified by a number of numerical tests using reference solutions.

\section{Theory}

We utilize the correspondence principle for electromagnetic wave and diffusion fields
\citep{lee1989new,de1996general,Mittet_2010_HFD} to calculate the CSEM response efficiently.
The key to a high-order local truncation error on a nonuniform grid relies on the solutions to a
Vandermonde system giving derivative operator coefficients that adapts to local grid properties. 
The stability condition is
then established for FDTD modelling on the nonuniform grid.

The Maxwell equations in a quasi-static regime (i.e., with negligible effect of displacement currents)
are written in the time domain as
\begin{eqnarray}\label{eq:time}
  \nabla \times \mathbf{E} +\mu\partial_t \mathbf{H} &=& -\mathbf{M}, \\ \nonumber
  -\nabla \times \mathbf{H}  +\bsig \mathbf{E} &=& -\mathbf{J},
\end{eqnarray}
or in the frequency domain as
\begin{eqnarray}\label{eq:freq}
  \nabla \times \mathbf{E} -\mathrm{i}\omega\mu \mathbf{H} &=& -\mathbf{M}, \\ \nonumber
  -\nabla \times \mathbf{H} +\bsig \mathbf{E} &=& -\mathbf{J},
\end{eqnarray}
where $\mathbf{E}=(E_x, E_y, E_z)^T$ and $\mathbf{H}=(H_x, H_y, H_z)^T$ are electric and magnetic
fields. The magnetic permeability is $\mu$.   
The conductivity
is a symmetric $3\times 3$ tensor: $\sigma_{ij}=\sigma_{ji},\; i,j\in\{x,y,z\}$. An isotropic medium
means that only the diagonal elements of the conductivity
tensor are non-zeros and the same in all directions: $\sigma_{xx}=\sigma_{yy}=\sigma_{zz}$; $\sigma_{ij}=0, i\neq j $.
The vertical transverse isotropic (VTI) medium implemented here still has only diagonal 
elements, but the vertical and the horizontal conductivities may differ, i.e.,
$ \sigma_h:=\sigma_{xx}=\sigma_{yy}$, $\sigma_v=\sigma_{zz}$. 
We use the following Fourier transform convention, $\partial_t \leftrightarrow -\mathrm{i}\omega$.

To speed up the FDTD modelling, we transform the above system from the diffusion to the wave domain,
following \citet{Mittet_2010_HFD}. The idea is to define a fictitious dielectric permittivity
in equation \ref{eq:freq} as $\sigma = 2\omega_0 \varepsilon$, yielding 
\begin{eqnarray}
  \nabla\times \mathbf{E} -\mathrm{i}\omega\mu \mathbf{H} &=& -\mathbf{M}, \\ \nonumber
  -\nabla\times \mathbf{H} +2\omega_0\varepsilon \mathbf{E} &=& -\mathbf{J},
\end{eqnarray}
which gives the following relation after multiplying the second equation with $\sqrt{-\mathrm{i}\omega/2\omega_0}$
\begin{eqnarray}\label{eq:mittet}
  \nabla\times \mathbf{E} +\underbrace{\sqrt{-\mathrm{i} 2\omega \omega_0}}_{-\mathrm{i}\omega'}\mu \underbrace{\sqrt{\frac{-\mathrm{i}\omega}{2\omega_0}} \mathbf{H}}_{\mathbf{H}'} &=& -\mathbf{M}, \\ \nonumber
  -\nabla\times\underbrace{ \sqrt{\frac{-\mathrm{i}\omega}{2\omega_0}} \mathbf{H}}_{\mathbf{H}'} +\underbrace{\sqrt{-\mathrm{i}2\omega\omega_0}}_{-\mathrm{i}\omega'}\varepsilon  \mathbf{E} &=& -\underbrace{\sqrt{\frac{-\mathrm{i}\omega}{2\omega_0}} \mathbf{J},}_{\mathbf{J}'}
\end{eqnarray}
which translates into the wave and simulation domain as the time dependent system
\begin{eqnarray}
  \nabla\times \mathbf{E}' + \mu\partial_t \mathbf{H}' &=& -\mathbf{M}', \\ \nonumber
  -\nabla\times \mathbf{H}' + \varepsilon\partial_t \mathbf{E}' &=& -\mathbf{J}'.
\end{eqnarray}
We have introduced a prime to identify the fields in the wave domain.

From the electromagnetic fields in the wave domain, the frequency-domain fields can be computed on
the fly during modelling using the fictitious wave transformation, exemplified by the electric field here,
\begin{eqnarray}\label{eq:dtft}
  \mathbf{E}'(\mathbf{x},\omega') = \int_0^{T_{\max}} \mathbf{E}'(\mathbf{x},t)e^{\mathrm{i}\omega't} \mathrm{d}t,
\end{eqnarray}
where $T_{\max}$ is the final time until the field $E'(\mathbf{x},\omega')$ reaches its steady state 
and where
\begin{eqnarray}\label{eq:omega_}
  \omega' = (1+\mathrm{i})\sqrt{\omega\omega_0}.
\end{eqnarray}

In order to have results valid for the frequency domain we need to calculate the Green's functions.
We need the following relation,
\begin{eqnarray}\label{eq:correspondence}
  \mathbf{E}' &=& \mathbf{E}, \\ \nonumber
  \mathbf{H}' &=& \sqrt{\frac{-\mathrm{i}\omega}{2\omega_0}} \mathbf{H},\\ \nonumber
  \mathbf{M}' &=& \mathbf{M},\\ \nonumber
  \mathbf{J}' &= & \sqrt{\frac{-\mathrm{i}\omega}{2\omega_0}} \mathbf{J} .
\end{eqnarray}
The Green's functions are then obtained by normalizing the transformed electric and magnetic fields with the source current,
\begin{eqnarray}
  G_{kj}^{E|J}(\mathbf{x},\omega|\mathbf{x}_s) = \frac{E_k(\mathbf{x},\omega|\mathbf{x}_s)}{J_j(\omega)}
  = \sqrt{\frac{-\mathrm{i}\omega}{2\omega_0}}\frac{E_k'(\mathbf{x},\omega|\mathbf{x}_s)}{J_j'(\omega)} , \\ \nonumber
  G_{kj}^{H|J}(\mathbf{x},\omega|\mathbf{x}_s) = \frac{H_k(\mathbf{x},\omega|\mathbf{x}_s)}{J_j(\omega)}
  = \frac{H_k'(\mathbf{x},\omega|\mathbf{x}_s)}{J_j'(\omega)},
\end{eqnarray}
where $G_{kj}^{E|J}(\mathbf{x},\omega|\mathbf{x}_s)$ and $G_{kj}^{H|J}(\mathbf{x},\omega|\mathbf{x}_s)$ stand for the electrical and
magnetic Green's function for angular frequency $\omega$ at spatial location $\mathbf{x}$ with the source located at $\mathbf{x}_s$.

Equation \ref{eq:mittet} is a pure wave-domain equation and the time integration can easily 
be discretized using the leap-frog method. We let the time be $t_n = n \Delta t$ with $n$ the
integer time variable and $\Delta t$ the time step. We also introduce $N = n + \frac{1}{2}$ such that,
\begin{eqnarray}\label{eq:timex}
  {\bHp}^{N} &=& {\bHp}^{N-1} + \Delta t\mu^{-1} (- \nabla\times {\bEp}^{n} - {\bMp}^{n}), \\ \nonumber
  {\bEp}^{n+1} &=& {\bEp}^{n} + \Delta t \epsilon^{-1} ( \nabla\times {\bHp}^{N} - {\bJp}^{N}).
\end{eqnarray}
The time integration of these equations is second-order accurate. It is shown in \citet{Mittet_2010_HFD}
that the calculation of the desired fields in the ``real world'' diffusive domain is independent
of the frequency content of the source term used for calculating the fictitious fields. We exploit
this fact and achieve good accuracy for the time integration by transmitting a low-frequency
signal in the fictitious wave domain. Here we are concerned with the spatial part of the
simulation scheme so we turn to this topic next.

We use a similar notation for the space variables as for the time variables where we write
$x_i = x_{i-1} +\Delta x_i$ where $\Delta x_i$ is the node separation between node $x_{i-1}$ and node $x_i$.
Likewise, we assume a forward staggered grid such that $x_I = x_{I-1}+ \Delta x_I$ 
where $\Delta x_I$ is the node separation between node $x_{I-1}$ and node $x_I$. 
For a uniform staggered grid we have a constant node separation such that $\Delta x_I = \Delta x_i = \Delta x$ and
$I = i + \frac{1}{2}$.
The $y$ and $z$ directions can be described in the same way with lower case and upper case integer arguments.

Calculation of the partial derivative of the field $f(x)$ can, in the continuous case, be formulated as an
integral operator by
\begin{eqnarray}\label{eq:alf}
  \partial_x f(x) = \int_{-\infty}^{\infty} \mathrm{d}x' f(x+x') \{ -\partial_{x'} \delta(x') \}
  = \int_{-\infty}^{\infty} \mathrm{d}x' f(x+x') \alpha(x').
\end{eqnarray}
The discrete formulation, with $f(x_i) := f(i)$, is,
\begin{eqnarray}\label{eq:dalf}
  \partial_x f(x_i) \approx D_x f(i)= \sum_{l=-L}^L f(i+l) \alpha_l(i) ,
\end{eqnarray}
where $\alpha_l(i)$ is a band-limited approximation to the operator $\alpha(x')$ in equation \ref{eq:alf}.
The half length of the operator is $L$. 
The argument $i$ is used to explicitly show that this operator will vary with location for a
nonuniform grid.
For the staggered grid, we can then define discretized forward, $ D_x^+ $, and backward,
$ D_x^- $, derivative operators as,
\begin{eqnarray}\label{eq:gandalf}
  D_x^+ f(i) = \partial_x f(I) = \sum_{l=1}^L f(i+l) \alpha_l(i) - f(i-l+1) \alpha_{-l}(i), \\ \nonumber
  D_x^- f(I) = \partial_x f(i)  = \sum_{l=1}^L f(I+l-1) \alpha_l(I) - f(I-l) \alpha_{-l}(I), 
\end{eqnarray}
which is the form we need for nonuniform grids and which is investigated here. The $f^{'}$ implies spatial
derivative in equation \ref{eq:gandalf}.

The operator simplifies for uniform grids where the operator becomes independent of spatial location
such that $ \alpha_l(i) = \alpha_{-l}(i) = \alpha_l$ 
\begin{eqnarray}\label{eq:simplyalf}
  D_x^+ f(i) = \partial_x f(I) = \sum_{l=1}^L ( f(i+l) - f(i-l+1) ) \alpha_{l}, \\ \nonumber
  D_x^- f(I) = \partial_x f(i) = \sum_{l=1}^L ( f(I+l-1) -f(I-l)  ) \alpha_{l}. 
\end{eqnarray}
If we use $L=1$ we have that $ \alpha_1 = 1/\Delta x $ and equation \ref{eq:simplyalf}
formulates the well known second-order accurate partial derivative operations,   
\begin{eqnarray}
  D_x^+ f(i) = \partial_x f(I)  = ( f(i+1) - f(i) ) / \Delta x, \\ \nonumber
  D_x^- f(I) = \partial_x f(i)  = ( f(I) - f(I-1) ) / \Delta x.
\end{eqnarray}

Let us assume a non-magnetic subsurface so that $\mu$ has the same value as in the vacuum, $\mathbf{M}=0$, while $\epsilon$ is a
diagonal tensor $\epsilon = \mbox{diag}(\epsilon_{ii}), i=x,y,z$. The staggering is as in \citet{Mittet_2010_HFD},
\begin{eqnarray}
  & & {\Hp}_x^{N}(i,J,K), \quad   {\Hp}_y^{N}(I,j,K), \quad     {\Hp}_z^{N}(I,J,k), \\ \nonumber
  & & {\Ep}_x^{n}(I,j,k), \quad    {\Ep}_y^{n}(i,J,k), \quad   {\Ep}_z^{n}(i,j,K), \\ \nonumber
  & & {\Jp}_x^{n}(I,j,k), \quad    {\Jp}_y^{n}(i,J,k), \quad   {\Jp}_z^{n}(i,j,K), \\ \nonumber
  & & \veps_{xx}(I,j,k), \quad     \veps_{yy}(i,J,k), \quad  \veps_{zz}(i,j,K), 
\end{eqnarray}
and the scheme implemented is, 
\begin{eqnarray}
  {\Hp}_x^{N}  &=& {\Hp}_x^{N-1} - \frac{\Delta t}{\mu}  ( D_y^+ {\Ep}_z^n  - D_z^+ {\Ep}_y^n ), \\ \nonumber
  {\Hp}_y^{N}  &=& {\Hp}_y^{N-1} - \frac{\Delta t}{\mu}  ( D_z^+ {\Ep}_x^n  - D_x^+ {\Ep}_z^n ), \\ \nonumber
  {\Hp}_z^{N}  &=& {\Hp}_z^{N-1} - \frac{\Delta t}{\mu}  ( D_x^+ {\Ep}_y^n  - D_y^+ {\Ep}_x^n ), \\ \nonumber
  {\Ep}_x^{n+1} &=& {\Ep}_x^{n} + \frac{\Delta t}{\epsilon_{xx}} ( D_y^- {\Hp}_z^N - D_z^- {\Hp}_y^N- {\Jp}_x^N ), \\ \nonumber
  {\Ep}_y^{n+1} &=& {\Ep}_y^{n} + \frac{\Delta t}{\epsilon_{yy}} ( D_z^- {\Hp}_x^N - D_x^- {\Hp}_z^N- {\Jp}_y^N ), \\ \nonumber
  {\Ep}_z^{n+1} &=& {\Ep}_z^{n} + \frac{\Delta t}{\epsilon_{zz}} ( D_x^- {\Hp}_y^N - D_y^- {\Hp}_x^N- {\Jp}_z^N ).
\end{eqnarray}

The conventional second-order staggered grid FDTD scheme discretizes the spatial derivatives as
follows \citep{newman1995frequency,mulder2006multigrid}:
\begin{eqnarray}
  & &   D_y^- H'_z = \frac{H'_z(I,J,k) - H'_z(I,J-1,k) }{\Delta y_{J}} , \quad 
  D_z^- H'_y = \frac{H'_y(I,j,K) - H'_y(I,j,K-1) }{\Delta z_{K}},  \nonumber \\
  & &  D_z^- H'_x = \frac{H'_x(i,J,K) - H'_x(i,J,K-1) }{\Delta z_{K}} , \quad  
  D_x^- H'_z = \frac{H'_z(I,J,k) - H'_z(I-1,J,k) }{\Delta x_{I}},   \nonumber \\
  & &  D_x^- H'_y = \frac{H'_y(I,j,K) - H'_y(I-1,j,K) }{\Delta x_{I}}, \quad 
  D_y^- H'_x = \frac{H'_x(i,J,K) - H'_x(i,J-1,K) }{\Delta y_{J}},    \nonumber \\
  & &   D_y^+  E'_z = \frac{E'_z(i,j+1,K) - E'_z(i,j,K) }{\Delta y_j}, \quad 
  D_z^+  E'_y = \frac{E'_y(i,J,k+1) - E'_y(i,J,k) }{\Delta z_k},   \nonumber \\
  & &  D_z^+  E'_x = \frac{E'_x(I,j,k+1) - E'_x(I,j,k) }{\Delta z_k}, \quad 
  D_x^+  E'_z = \frac{E'_z(i+1,j,K) - E'_z(i,j,K) }{\Delta x_i},   \nonumber \\
  & &  D_x^+  E'_y = \frac{E'_y(i+1,J,k) - E'_y(i,J,k) }{\Delta x_i}, \quad 
  D_y^+  E'_x = \frac{E'_x(I,j+1,k) - E'_x(I,j,k) }{\Delta y_j},
\end{eqnarray}
where $\Delta x_i$, $\Delta y_j$ and $\Delta z_k$ are distances between nodes on the reference grid,
while $\Delta x_{I}$, $\Delta y_{J}$ and $\Delta z_{K}$ are the distance
between the grid points $(I,j,k)$ and  $(I-1,j,k)$, the distance between the grid points
$(i,J,k)$ and  $(i,J-1,k)$, and the distance between the grid points $(i,j,K)$ and  $(i,j,K-1)$.
This scheme is second-order accurate on a uniform grid with $\Delta x_i = \Delta x_{I}$ and likewise
for the other spatial directions. The scheme is consistent with using equation \ref{eq:simplyalf}
with $L=1$ to approximate the derivatives.

A standard discretization method for nonuniform grids is to use the same formulation as above, but
where the node distance may vary along the same spatial direction. It is well known that there
are accuracy issues with this implementation. 
To achieve second-order accuracy using
equation \ref{eq:simplyalf} is not possible. The midpoints between nodes do not align after
going from the reference grid to the staggered grid and back again.
Since $\Delta x_i \neq \Delta x_{I}$, we have a situation where the cell center does not
match on a nonuniform grid, as is illustrated in Figure~\ref{fig:figure1}b. 
Consequently, the local truncation error of the resulting scheme can only reach first order.
Reduced accuracy will also be a problem if we implement a high-order FDTD scheme ($L>1 $) on a nonuniform grid,
using derivative-operator coefficients designed for a regular grid. 
This is unfortunate since the nonuniform grid is potentially attractive for efficient modelling due to
significant reduction of the number of grid point. However, good accuracy can be restored if equation \ref{eq:gandalf}  
is used instead of equation \ref{eq:simplyalf}. The problem that remains is to calculate the
operator coefficients for equation \ref{eq:gandalf}. 

\subsection{Vandermonde matrix}

\begin{figure}
  \centering
  \includegraphics[width=0.5\textwidth]{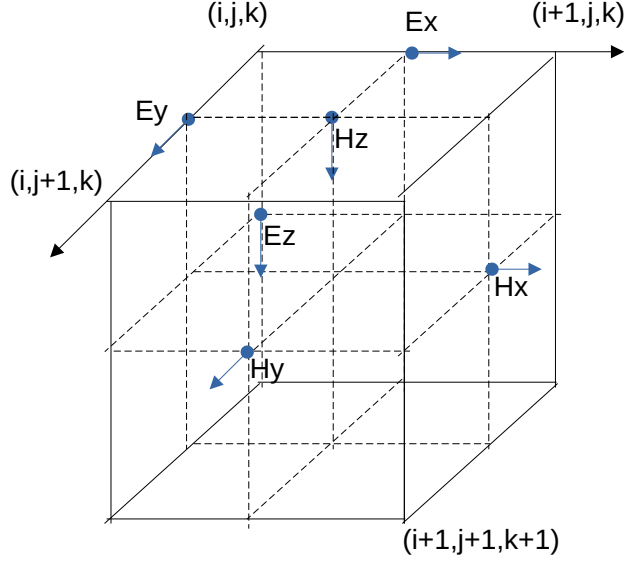}
  \caption{Electrical and magnetic fields on staggered grid \citep{Yee_1966_NSI}.}\label{fig:figure2}
\end{figure}

The major difference between FDTD implementations on a uniform grid and on a nonuniform grid lies in
the design of the spatial-derivative operator coefficients.
The position of each field component on the nonuniform staggered grid has been illustrated in
Figure~\ref{fig:figure2}, which is similar to the staggered FDTD on a uniform grid. 

In order to compute the electromagnetic field as well as its derivatives with arbitrary
grid spacing, we have to do a polynomial interpolation using a number of knots
$x_0, x_1,\cdots, x_n$. According to the Taylor expansion, we have
\begin{eqnarray}
  f(x_i) &=& f(x) + f'(x)(x_i-x) + \frac{1}{2}f''(x)(x_i-x)^2 + \cdots +
  \frac{1}{n!}f^{(n)}(x)(x_i-x)^n+\cdots \\ \nonumber
  i &=& 0,1,\cdots,n.
\end{eqnarray}
By defining $a_i(x):=f^{(i)}(x)/i!$, we end up with a polynomial of the Newton form,
\begin{eqnarray}
  f(x_i) &=& a_0(x) + a_1(x)(x_i-x) + a_2(x)(x_i-x)^2 + \cdots +a_n(x)(x_i-x)^n+\cdots \\ \nonumber
  i &=& 0,1,\cdots,n. 
\end{eqnarray}
Let us consider $n+1$ distinct nodes $x_0,x_1,\cdots, x_n$ and drop the terms
$O((x_i-x)^{n+1})$. This builds a matrix system
\begin{eqnarray}
  \underbrace{  \begin{pmatrix}
      f(x_0)\\
      f(x_1)\\
      \vdots\\
      f(x_n)
  \end{pmatrix}}_{\mathbf{f}}
  =\underbrace{\begin{pmatrix}
      1 & x_0-x & (x_0-x)^2 & \cdots & (x_0-x)^n\\
      1 & x_1-x & (x_1-x)^2 & \cdots & (x_1-x)^n\\
      \cdots\\
      1 & x_n-x & (x_n-x)^2 & \cdots & (x_n-x)^n\\
  \end{pmatrix}}_{\mathbf{V}^T(x_0-x,\cdots,x_n-x)} \underbrace{\begin{pmatrix}
      a_0(x)\\
      a_1(x)\\
      \vdots\\
      a_n(x)
  \end{pmatrix}}_{\mathbf{a}},
\end{eqnarray}
where $\mathbf{V}^T(x_0-x,\cdots,x_n-x)$ is the transpose of a Vandermonde matrix
determined by $x_0-x,\cdots,x_n-x$. The above expression implies that the function $f(x)$
and its derivatives up to the $n$-th order at arbitrary location $x$ can be found by inverting
the Vandermonde matrix: $(f(x), f'(x),\cdots, f^{(n)}(x)/n!)^T = \mathbf{a}=[\mathbf{V}^T]^{-1}\mathbf{f}$.

It is well known that the Vandermonde matrix is highly ill-conditioned and direct matrix
inversion by Gaussian elimination should be avoided due to numerical instabilities when the matrix
size becomes large. Fortunately, there exists an efficient algorithm based on the method of
\citet{bjorck1970solution}  to invert the Vandermonde matrix. In fact, the algorithm circumvents
the curse of severe ill-conditioning of the Vandermonde matrix to arrive at arbitrarily high
accuracy for the inversion \citep{demmel2005accurate}.
Compared with Gauss-elimination of complexity $O(n^3)$ , the Vandermonde matrix inversion
algorithm reduces the computational complexity to $O(n^2)$. 
Because the elements of the Vandermonde matrix are fully determined by the interpolation nodes,
there is no need to explicitly construct the matrix and store it before inversion.
The detailed implementation of  this fast algorithm is available in
\citet[Algorithm 4.6.1]{Golub_1996_MATCOMP}.

Let the $i$-th row, $j$-th column of the inverse matrix $[\mathbf{V}^T]^{-1}$
be $w_{ij}$, i.e., $([\mathbf{V}^T]^{-1})_{ij}=w_{ij}, i,j=0,\cdots,n$. It also follows that
\begin{eqnarray}\label{eq:amatrix}
  \underbrace{  \begin{pmatrix}
      a_0(x)\\
      a_1(x)\\
      \vdots\\
      a_n(x)
  \end{pmatrix}}_{\mathbf{a}}:=\underbrace{\begin{pmatrix}
      w_{00} & w_{01} & \cdots w_{0n}\\
      w_{10} & w_{11} & \cdots w_{1n}\\
      \vdots\\
      w_{n0} & w_{n1} & \cdots w_{nn}\\
  \end{pmatrix}}_{[\mathbf{V}^T]^{-1}}\underbrace{\begin{pmatrix}
      f(x_0)\\
      f(x_1)\\
      \vdots\\
      f(x_n)
  \end{pmatrix}}_{\mathbf{f}}.
\end{eqnarray}
The $i$-th row gives the explicit expression to find the $i$-th derivative 
\begin{eqnarray}\label{eq:interp}
  \frac{1}{i!} f^{(i)}(x) = a_i(x)= \sum_{j=0}^n  w_{ij} f(x_j)= \langle w_{i\cdot}|\mathbf{f}\rangle, \quad i=0,\cdots,n.
\end{eqnarray}

The first row of the matrix $[\mathbf{V}^T]^{-1}$ (i.e., $w_{0j}$, $j=0,\cdots, n$) is $f(x)$
which in analogy with equation \ref{eq:alf} can be written in discrete form as
\begin{eqnarray}\label{eq:bet}
  f(x) = \int_{-\infty}^{\infty} \mathrm{d}x' f(x+x') \{ \delta(x') \}
  = \int_{-\infty}^{\infty} \mathrm{d}x' f(x+x') \beta(x').
\end{eqnarray}
The discrete formulation is,
\begin{eqnarray}\label{eq:dbet}
  f(x_i) = \sum_{l=-L}^L f(i+l) \beta_l(i) ,
\end{eqnarray}
where $\beta_l(i)$ are the coefficients for an interpolation operator adapted to a nonuniform grid
and identical to the  $w_{0j}$ coefficients.
The second row of the matrix $[\mathbf{V}^T]^{-1}$ (i.e., $w_{1j}$, $j=0,\cdots, n$) is $\partial_x f(x)$
and the continuous and discrete representations are  given in equations \ref{eq:alf} and \ref{eq:dalf}.
The $\alpha_l(i)$ coefficients are for a derivative operator adapted to a nonuniform grid
and identical to the $w_{1j}$ coefficients. To be explicit, 
consider the staggered finite-difference approximation of the first derivatives in $x$
direction using $2L$ non-equidistant nodes.  The finite-difference coefficients
$\alpha_l(x_{i})$ and $\alpha_l(x_{I})$,  $l=-L+1,\cdots, L$ are the 2nd row of the inverse of the
matrices $\mathbf{V}^T(x_{i+L}-x_I,\cdots,x_{i-L+1}-x_I)$ and $\mathbf{V}^T(x_{I+L}-x_i,\cdots,x_{I-L+1}-x_i)$.
Using the $2L$ nodes, we achieve accuracy up to $2L$-th order in space.

In general we find that the operator coefficients
(interpolation weights) for $f^{(i)}(x)$ are $i!w_{ij}$.
Given the points $x_0,\cdots,x_n$ and $x$, the Vandermonde matrix is determined and the
operator coefficients can be calculated. For the simulation we only need the derivative-operator coefficients. 
It is noteworthy to mention that the Vandermonde matrix must be non-singular to be inverted.
For the derivative-operator coefficients we have that the Vandermonde matrix is non-singular
by construction of the staggered grids. A finite-difference scheme
on a staggered grid implies that the node $x_I$, whose derivative is computed, will stay in the middle
between the selected nodes $x_i$, which ensures that the resulting Vandermonde matrix is invertible.

The interpolation operators are useful for recording fields at arbitrary locations and for the
distribution of source contributions \citep{mittet2017internal}. Just as for the derivative-operator coefficients,
the interpolation coefficients can be pre-calculated and then reused every time step.
For the interpolation operator we may find that some of the coordinates $x_i$ matches the interpolation points
$x$, the inversion of Vandermonde matrix is then not necessary.
This might happen when the source or receiver positions coincide with a finite-difference node.
In this case, the interpolation weights $w_{0i}$ to evaluate $f(x)$ should be exactly 1 at $x_i$
and $0$ elsewhere. In order to do 2D/3D simulation on nonuniform grid, multidimensional
interpolant is simply constructed by tensor products of many 1D interpolants.

The above procedure is significant as it allows us to use arbitrarily high-order
finite-difference scheme to accurately compute the electromagnetic fields and their derivatives,
typically with arbitrary grid spacing in the rectilinear grid. This opens the door for CSEM
modelling using high-order FDTD on a nonuniform grid in a consistent framework.
The computed finite-difference coefficients may also be used to do high-order frequency-domain
modelling on a nonuniform grid, while the resulting sparse banded matrix has to be solved
accurately if sufficient computational resources are available.

\subsection{Stability condition}

Let us write down the FDTD scheme in equation \ref{eq:timex} without source terms as follows:
\begin{eqnarray}
  \begin{cases}
    \mathbf{E}^{'n+1} =\mathbf{E}^{'n} + \Delta t \epsilon^{-1} \nabla\times \mathbf{H}^{'N}\\
    \mathbf{H}^{'N+1} = \mathbf{H}^{'N} -\Delta t \mu^{-1}\nabla\times \mathbf{E}^{'n+1}
  \end{cases},
\end{eqnarray}
leading to
\begin{eqnarray}
  \begin{split}
    \begin{bmatrix}
      \mathbf{E}^{'n+1}\\
      \mathbf{H}^{'N+1}
    \end{bmatrix}   =\underbrace{\begin{bmatrix}
        \mathbf{I} & \Delta t \epsilon^{-1} \nabla\times\\
        -\Delta t \mu^{-1}\nabla\times & \mathbf{I}-\Delta t^2 \epsilon^{-1}\mu^{-1}\nabla\times\nabla\times
    \end{bmatrix}}_{\mathbf{A}}\begin{bmatrix}
      \mathbf{E}^{'n}\\
      \mathbf{H}^{'N}
    \end{bmatrix}.
  \end{split}
\end{eqnarray}
The numerical stability requires the eigenvalues of the amplification matrix $\mathbf{A}$
to be less than or equal to 1. Assume the eigenvalue decomposition for the amplification
matrix is $\mathbf{A} = \mathbf{\bar{V}}\mathbf{\Lambda} \mathbf{\bar{V}}^T$,
where $\mathbf{\bar{V}}=(\mathbf{V}_E,\mathbf{V}_H)^T$ is an unitary matrix
such that $\mathbf{\bar{V}}^T\mathbf{\bar{V}}=\mathbf{I}$.
Then we have $\mathbf{A}\mathbf{V}=\mathbf{V}\mathbf{\Lambda}$, yielding
\begin{eqnarray}
  \begin{bmatrix}
    \mathbf{I} & \Delta t \epsilon^{-1} \nabla\times\\
    -\Delta t \mu^{-1}\nabla\times & \mathbf{I}-\Delta t^2 \epsilon^{-1}\mu^{-1}\nabla\times\nabla\times
  \end{bmatrix}\begin{bmatrix}
    \mathbf{V}_E\\
    \mathbf{V}_H
  \end{bmatrix} = \begin{bmatrix}
    \mathbf{V}_E\\
    \mathbf{V}_H
  \end{bmatrix}\Lambda.
\end{eqnarray}
That is,
\begin{eqnarray}
  \begin{cases}
    \Delta t\epsilon^{-1}\mu^{-1}\nabla\times \mathbf{V}_H = \mathbf{V}_E(\mathbf{\Lambda}-\mathbf{I})\\
    -\Delta t \mu^{-1}\nabla\times \mathbf{V}_E -\Delta t^2\epsilon^{-1}\mu^{-1}\nabla\times\nabla\times \mathbf{V}_H = \mathbf{V}_H(\mathbf{\Lambda}-\mathbf{I}).
  \end{cases}
\end{eqnarray}
Multiplying the second sub equation $\mathbf{\Lambda}-\mathbf{I}$ from the right and inserting the first sub equation gives
\begin{equation}
  -\Delta t^2 \mu^{-1}\epsilon^{-1} \nabla\times\nabla\times \mathbf{V}_H \mathbf{\Lambda} = \mathbf{V}_H(\mathbf{\Lambda}-\mathbf{I})^2.
\end{equation}
Denote $\mathbf{V}_{H,j}$ the $j$th column of $\mathbf{V}_H$ and $\lambda_j$ the $j$th eigenvalue in $\mathbf{\Lambda}$. The above equation reads
\begin{equation}
  \mu^{-1}\epsilon^{-1} \nabla\times\nabla\times \mathbf{V}_{H,j}=-\frac{ (\lambda_j-1)^2}{\Delta t^2\lambda_j} \mathbf{V}_{H,j},
\end{equation}
which shows that $-\frac{ (\lambda_i-1)^2}{\Delta t^2\lambda_i} $ is the eigenvalue of the matrix $(\mu\epsilon)^{-1}\nabla\times\nabla\times$ associated with the eigenvector $V_{H,j}$. This leads to
\begin{equation}
  (\lambda_j^2 + (-2+\Delta t^2 c^2\nabla\times\nabla\times)\lambda_j + 1) \mathbf{V}_{H,j}
  =(\lambda_j^2 -(2+\Delta t^2 c^2\Delta)\lambda_j + 1) \mathbf{V}_{H,j}= 0,
\end{equation}
where we denote $c:=1/\sqrt{\mu\epsilon}$ and have applied $\nabla\times\nabla\times \mathbf{F} = \nabla\nabla\cdot \mathbf{F} - \nabla\cdot\nabla \mathbf{F}=-\Delta \mathbf{F}$ due to Gauss theorem $\nabla\cdot \mathbf{F} =0$, $\mathbf{F}=\mathbf{E},\mathbf{H}$ in the homogeneous, source free medium. The roots of the above equation are 
\begin{equation}
  \lambda_{j;1,2} =  1+\frac{\Delta t^2c^2\Delta }{2} \pm \frac{\mathrm{i}}{2}\sqrt{-\Delta t^2 c^2\Delta (4+\Delta t^2c^2\Delta )},
\end{equation}
which requires the following condition to be satisfied
\begin{equation}
  0\leq -\Delta t^2 c^2 \Delta  \leq 4,
\end{equation}
in order to ensure  $|\lambda_{j;1,2}|\leq 1$. Finally, we arrive at the same stability condition as
equation 41 of \citet{Mittet_2010_HFD},
\begin{equation}
  \Delta t c_{\max}\sqrt{(D_x^{\max})^2 + (D_y^{\max})^2 + (D_z^{\max})^2} \leq 2,
\end{equation}
where $D_x^{\max}$, $D_y^{\max}$ and $D_z^{\max}$ are the maximum value of the the discretized first derivatives along $x$, $y$ and $z$ directions. Let us emphasize this condition applies to both uniform and nonuniform grid. The difference lies in the spatial derivative operator.

To proceed with the stability analysis, we represent the fields on the grid via time harmonic plane waves
\begin{equation}
  \mathbf{E}',\mathbf{H}' \propto e^{-\mathrm{i}(\omega t - k_x  x - k_y  y - k_z z)},
\end{equation}
where the amplitude has been omitted. Equation \ref{eq:gandalf} becomes
\begin{equation}
  \begin{cases}
    D_x^+ u (x_I) = (\alpha_L(x_I) e^{\mathrm{i}k_x (x_{i+L}-x_I)} + \alpha_{L-1}(x_I) e^{\mathrm{i}k_x (x_{i+L-1}-x_I)} + \cdots + \alpha_{-L+1}(x_I)e^{\mathrm{i}k_x (x_{i-L+1}-x_I)} )u(x_I)\\
    D_x^- u(x_i) = (\alpha_L(x_i) e^{\mathrm{i}k_x (x_{I+L-1}-x_i)} + \alpha_{L-1}(x_i) e^{\mathrm{i}k_x (x_{I+L-1}-x_i)} + \cdots + \alpha_{-L+1}(x_i)e^{\mathrm{i}k_x (x_{I-L}-x_i)} )u(x_i)
  \end{cases}.
\end{equation}
Hence, we end up with the maximum possible values for discrete first derivative operators
\begin{equation}
  D_x^{\max}  = \max\left(\sum_{l=-L+1}^L |\alpha_l(x_I)|,\sum_{l=-L+1}^L |\alpha_i(x_l)|\right)
\end{equation}
and similar estimations for $D_y^{\max}$ and $D_z^{\max}$. In case of a uniform grid, the above expressions becomes much simpler
\begin{displaymath}
  \begin{cases}
    D_x^+ u(x_I) = (\alpha_L e^{\mathrm{i}k_x (L-1/2)\Delta x} + \alpha_{L-1} e^{\mathrm{i}k_x (L-3/2)\Delta x} + \cdots + \alpha_{-L+1}e^{\mathrm{i}k_x (-L+1/2)\Delta x}u(x_I)\\
    D_x^- u(x_i) =  (\alpha_L e^{\mathrm{i}k_x (L-1/2)\Delta x} + \alpha_{L-1} e^{\mathrm{i}k_x (L-3/2)\Delta x} + \cdots + \alpha_{-L+1}e^{\mathrm{i}k_x (-L+1/2)\Delta x}u(x_i)\\
  \end{cases},
\end{displaymath}
where $\Delta x$ stands for the uniform grid spacing in $x$ direction, while the coefficients $\alpha_l$ (which can be computed by inverting a Vandermonde matrix system  according to Appendix A\ref{appendix:fdcoeff}) are independent of the location $x_i$.

\subsection{Grid stretching}

Our finite-difference modelling is carried out on a rectilinear mesh, which can be generated
from the tensor (outer) product of 1D non-equispaced meshes.
We use the geometrical progression to generate the 1D nonuniform grid, following the work of \citet[ Appendix C]{mulder2006multigrid}.
This is also often referred to as power law grid stretching since the cell sizes stretch
exponentially to guarantee a smooth extension of the grid.

Assume we have the total grid length $L_x$ divided into $n$ intervals ($n+1$ nodes) with a
common ratio $r>1$. Denote the smallest interval $\Delta x=x_1-x_0$. Thus, the relation
between $L_x$ and $\Delta x$ is
\begin{equation}\label{eq:optnugrid}
  L_x = (x_1-x_0)  + (x_2-x_1) + \cdots + (x_n-x_{n-1})
  = \Delta x(1 + r + \cdots + \cdot r^{n-1})= \Delta x \frac{r^n-1}{r-1}.
\end{equation}
Given the total distance $L_x$, the smallest grid spacing $\Delta x$ and the stretching factor $r$,
we can compute an approximate value for the number of nodes
$ n = \left\lceil \frac{  \ln(1 + \frac{L_x}{\Delta x}(r -1)) }{\ln(r)} \right\rceil$ following \citet{mulder2006multigrid},
where $\left\lceil\cdot\right\rceil$ takes the ceiling integer value. This strategy yields approximate
solution as the value of $L_x$ is not exactly preserved.

Due to the stability requirement and the resulting computational cost in the modelling,
we are restricted to the smallest interval $\Delta x$ and a given number of intervals $n$ to
discretize over a certain distance $L_x$. The question boils down to finding  an optimal growth factor $r$.
This problem is more complicated since equation \ref{eq:optnugrid} does not yield an explicit
expression for the stretching factor $r$. 

The relation in equation \ref{eq:optnugrid} is equivalent to 
\begin{equation}\label{eq:rr}
  r= \underbrace{\left(\frac{L_x}{\Delta x}(r-1) +1\right)^{\frac{1}{n}}}_{g(r)},
\end{equation}
which inspires us to carry out a number of fixed point iterations until convergence:
\begin{equation}\label{eq:fixed}
  r^{k+1} = g(r^k), \quad k=0,1,\cdots.
\end{equation}
Assume $r^*$ is the analytic solution such that $r^*=g(r^*)$. Thanks to Lagrange mean value
theorem, the error estimation at ($k+1$)-th iteration is linked with the error at $k$-th iteration via
\begin{equation}
  |e^{k+1}| = |r^{k+1}-r^*|=|g(r^k)-g(r^*)|=|g'(\xi)(r^k-r^*)|=|g'(\xi)| |e^k|, \quad \xi \mbox{ between } r^k \mbox{ and } r^*.
\end{equation}
It becomes evident that $e^k\rightarrow 0\; (k\rightarrow \infty)$ provided that $|g'(r)|<1$. 
Starting from any initial guess $r^0>1$, the fixed point iteration scheme in equation \ref{eq:fixed}
is guaranteed to converge since
\begin{displaymath}
  |g'(r)|= \frac{L_x}{n\Delta x}\left(\frac{L_x}{\Delta x}(r-1) +1\right)^{1/n-1}
  = \frac{1}{n}(r^{-(n-1)}+r^{-(n-2)}+ \cdots +1) <1,
\end{displaymath}
thanks to the relations in equations \ref{eq:optnugrid} and \ref{eq:rr}.

We note that within areas of constant $r$, the derivative-operator coefficients in equation \ref{eq:gandalf}
can be calculated from the $2L$ coefficients $a_{-l}$ and $a_l$, such that
\begin{eqnarray}
  \alpha_{-l}(i) &=& \frac{a_{-l}}{\Delta x \: r^i }, \\ \nonumber
  \alpha_{l}(i) &=& \frac{a_{l}}{\Delta x \: r^i }.  \\ \nonumber
\end{eqnarray}
Two solutions of the Vandermonde system are required in this case, one for operators valid on the reference grid and one for
operators valid on the staggered grid.

\subsection{Implementation}

A number of techniques have been applied to achieve efficient and accurate 3D CSEM simulation.
According to equation \ref{eq:dtft}, the frequency-domain CSEM response is obtained by a time
to frequency transform. The transform is using the complex frequency given in equation \ref{eq:omega_}.
This gives exponential damping of late arrivals as is discussed in more detail in \citet{mittet2015seismic}.
The lowest frequency experiences the least damping in this transform and by that requires the longest
simulation time. This allows us to bound the number of times steps to terminate the simulation
when the lowest frequency component has converged. The convergence means that the
frequency-domain field obtained by time integration has reached its steady state and later arrivals are
damped to the degree that they do not contribute to the time integral.

The source and receiver locations may be arbitrarily distributed over the whole computational
domain, not necessarily located at the nodes of the finite-difference grid. In case they do,
interpolation is not required and we directly take the field from the grid; otherwise,
we need interpolation operators extracted from the first row of the relevant
$[\mathbf{V}^T]^{-1}$ matrices.

The CFL condition to achieve stable FDTD modelling dictates the timestep for a given spatial sampling.
The air-wave travels at extremely high speed which does not allow us to use the local finite-difference
stencil to simulate it. The air-water boundary condition is implemented in the Fourier-wavenumber domain
following the method proposed by \citet{oristaglio1984diffusion} and \citet{wang1993finite} and using the extension
to high-order schemes given in \citet{Mittet_2010_HFD}.

To mimic wave propagation in unbounded domain, we use a truncated computation mesh surrounded by
an artificial absorbing boundary using convolutional perfectly matched layer (PML) technique
\citep{roden2000convolution}, except the top air-water interface. To ease the implementation,
we extend the domain with equal spacing based on the last finite-difference cell size in the interior domain.
Our CPML implementation is very standard using all the parameter setups given in \citet{Komatitsch_2007_GEO}.

\section{Numerical examples}

We now present three examples to demonstrate the merits of high-order finite differences using nonuniform staggered grids. The first two examples are using 1D resistivity models under shallow water and deep water
scenarios. In the 1D case, the semi-analytic solution can be computed in the frequency-wavenumber domain as a reference to benchmark our results.
The third example takes into account varying seafloor topography, in which
a reference solution can be computed using \verb|emg3d| software \citep{werthmuller2019emg3d}.

In all of our modelling, we use x-directed electrical dipole source. 
The $E_x$ component at three commonly used frequencies - 0.25 Hz, 0.75 Hz and 1.25 Hz, are modelled.
In total 12 layers of CPML are sufficient to achieve nearly perfect absorbing effect.
The computed electromagnetic fields are normalized by the source current.
The convergence check has been conducted every 100 time steps to avoid modelling after the
frequency-domain EM fields stop evolving.
We examine the amplitude error by inspecting the ratio of the modelled field to the reference solution, $|E_x^{FD}|/|E_x^{ref}|$, which should be close to unity if the modelling is precise.
The phase difference is computed  by $\angle E_x^{FD}-\angle E_x^{ref}$ in degrees.

\subsection{Moving to high-order schemes}

\begin{figure}
  \centering
  \includegraphics[width=0.5\linewidth]{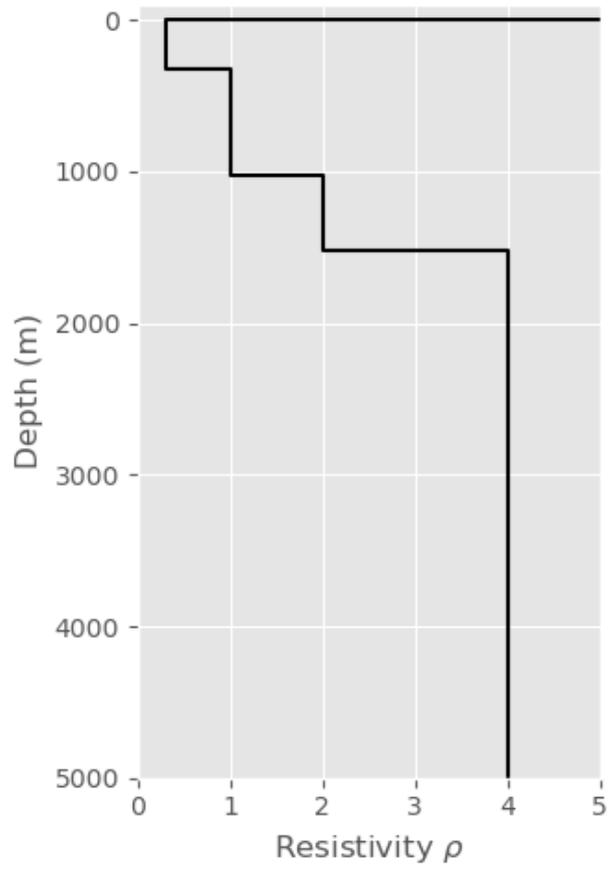}
  \caption{The resistivity model with air ($\rho=10^{12}$ $\Omega$m), shallow column (325 m) of 0.3 $\Omega$m
    and 3 sediment layers ($\rho=1$ $\Omega$m in [325, 1025] m; $\rho=2$ $\Omega$m in [1025, 1525] m; $\rho=4$ $\Omega$m downwards).}\label{fig:figure3}
\end{figure}

The model shown in Figure~\ref{fig:figure3} includes 5 layers: the top is the air,
then a 325 m layer of sea water with a resistivity of 0.3 $\Omega$m,
followed by three layers of formation with increasing resistivity in the depth direction. The whole model
extends down to 5 km depth below the sea surface.
A horizontal electrical dipole source is deployed at 275 m water depth, its lateral position
being in the middle of the model.
The EM fields are recorded using 201 receiver positions at the seafloor. The offsets range from -10 km to +10 km
(receiver separation equals 100 m). 

\begin{figure}
  \centering
  \includegraphics[width=\linewidth]{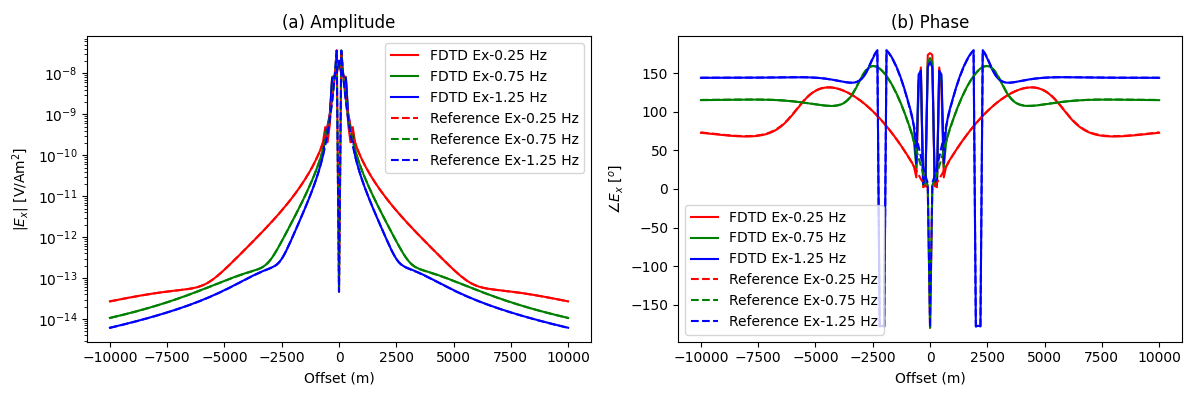}
  \caption{Comparison  between 2nd order FDTD (solid line) and reference solution (dash line) for 3D
    CSEM simulation in the shallow water scenario. The horizontal coordinates are offsets, while the vertical coordinates are (a) Amplitude; (b) Phase.}\label{fig:figure4}
\end{figure}

We first compare the 3D CSEM response simulated by FDTD
with the reference solution calculated using \verb|empymod| program \citep{werthmuller2017open}. To limit
the factors affecting the modelling accuracy, this experiment has been done using isotropic
modeling and uniform grid spacing ($\Delta x=\Delta y= 150$ m, $\Delta z=50$ m).
Panels (a) and (b) in Figure~\ref{fig:figure4} shows the amplitude and the phase of the modelled EM fields
by FDTD of 2nd order, which is highly consistent to the reference solution.
At the offset  beyond 1 km, the agreement between finite-difference and the analytical
solution is good. The kink in the transition of
the near-to-far offset in the  amplitude response for all frequencies is a manifestation 
of the strong air-wave effect.
In the very near field, where the receivers are located just below the transmitter,          
the phase exhibits a $180^o$ jump due to the change in the direction of the electric
field immediately below the electric dipole source.
Due to the extension of finite-difference stencil, this phase rollover becomes smeared
compared to the reference solution. For the same reason,  finite difference method introduces
significant errors in the near field in amplitude (which becomes more evident for
high-order schemes). Fortunately it is not an issue for practical 3D CSEM applications since
the CSEM inversion to deduce the subsurface resistivity is mainly driven by far offset
refractions. We thus focus on the analysis of the amplitude and phase error
beyond 1 km offset in the following.

\begin{figure}
  \centering
  \includegraphics[width=\linewidth]{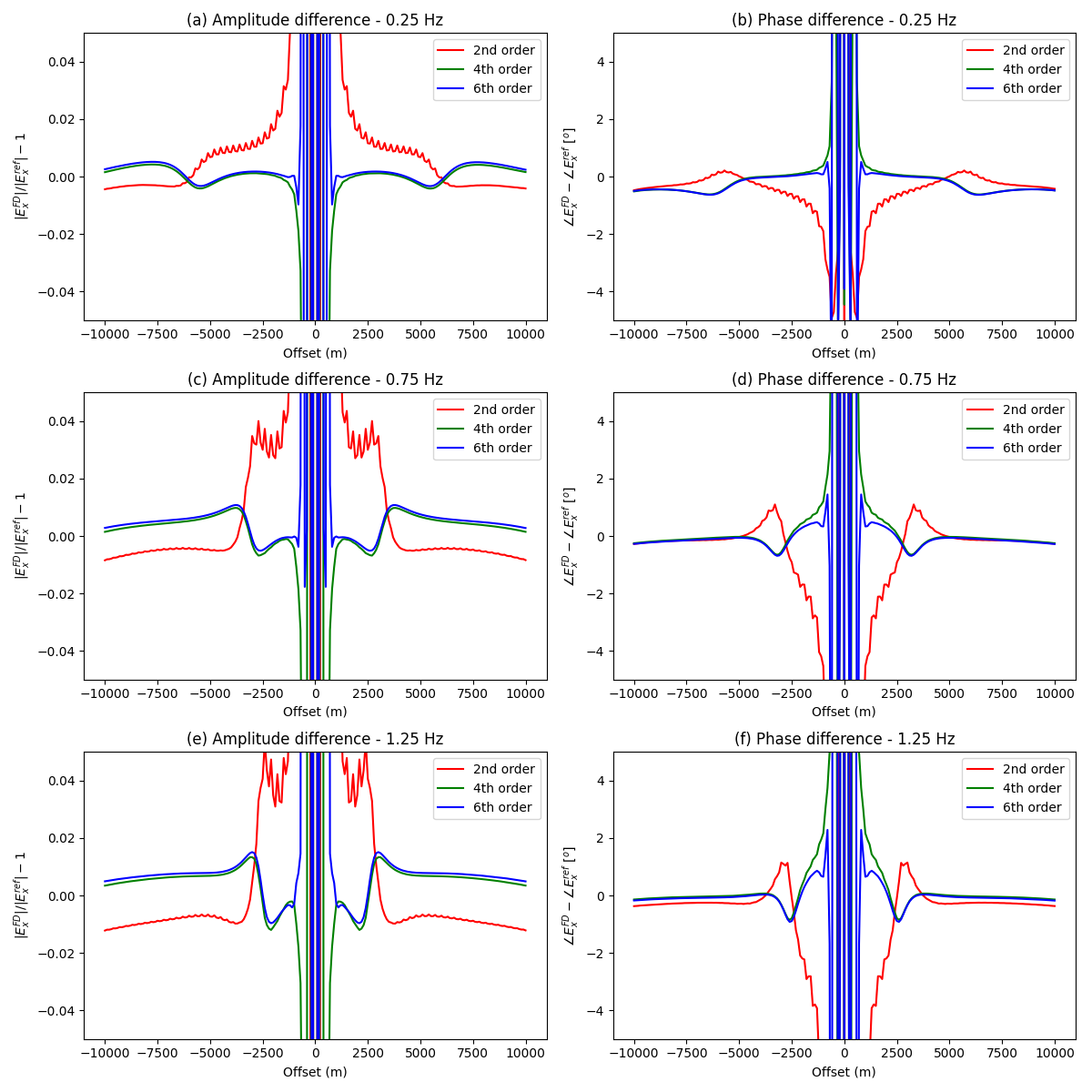}
  \caption{The amplitude and phase error of FDTD compared with 1D reference solution for 0.25 Hz (a,b), 0.75 Hz (c,d) and 1.25 Hz (e,f).}\label{fig:figure5}
\end{figure}

The modelling results by FDTD of higher orders are not displayed in Figure~\ref{fig:figure4},
since they are visually very similar to the reference solution. Instead, the amplitude and phase error are computed using 1D reference solution to examine the accuracy of our methods.
In Figure~\ref{fig:figure5}a, c and e, we clearly see that the 2nd order FDTD gives the largest amplitude error
for all frequencies; moving from 2nd order to 4th order significantly reduces the amplitude error, while moving to 6th order
behaves even better, although the accuracy improvement becomes less. The phase errors exhibit a similar behavior.

In principle, the computational cost of the 4th order and the 6th order finite differences will
be double and triple of that for the 2nd order scheme.  Table~\ref{table:time}
lists the CPU time for these modelling exercises, running on a laptop possessing Intel(R) Core(TM) i7-4710HQ CPU @ 2.50GHz.
It shows the increase of the
computing cost by increasing the order of FDTD scheme. The computing time of the 4th order
scheme is less than two times the computing time of the 2nd order scheme. Further increase of the FDTD order results
in significant increase of CPU time. Since moving to the 6th order scheme demands more
computation while the accuracy improvement is marginal, we stay with the 4th order
scheme from now on, as it gives the best compromise between increased accuracy and computational load.

\begin{table}
  \centering
  \caption{Comparison of computing time using FDTD of orders 2, 4 and 6.}\label{table:time}
  \begin{tabular}{l|c|c|c}
    \hline
    FDTD & Order-2 &  Order-4     &  Order-6\\
    \hline
    Time (s) &  422.3 & 715.5 & 1090.9\\
    \hline
  \end{tabular}
\end{table}

\subsection{The impact of nonuniform grid}

\begin{figure}
  \centering
  \includegraphics[width=0.85\textwidth]{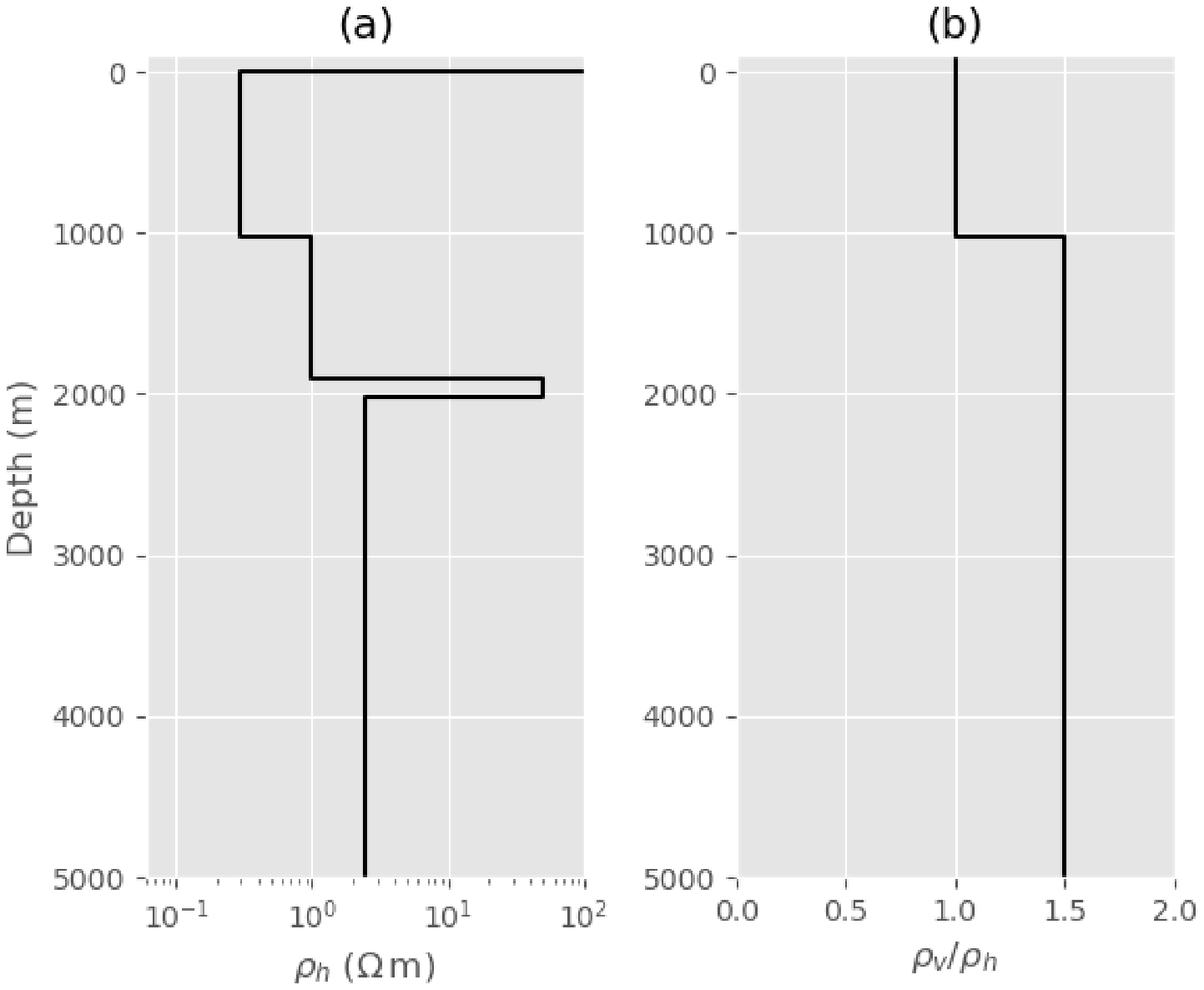}
  \caption{(a) The resistivity model in deep water: the 1020 m water column of 0.3 $\Omega$m followed by formation of $1 \Omega\mbox{m}$ down to 1900 m, and 120 m thickness of resistor of 50 $\Omega$m, while the background resistivity below the resistor is 2.5 $\Omega$m;    (b) VTI anisotropy below seabed is 1.5.}\label{fig:figure6}
\end{figure}

The above example demonstrates the importance of higher order scheme to achieve accurate CSEM modeling in the presence of strong air-wave. We now examine the impact of grid non uniformity in achieving
computational efficiency for high-order staggered grid FDTD. To get rid of the impact of air-wave, we consider a resistivity model in deep water scenario. As shown in Figure~\ref{fig:figure6}a, the model has 1020 m water column of 0.3 $\Omega$m, followed by formation of $1 \Omega\mbox{m}$ down to 1900 m, and 120 m thickness of resistor of 50 $\Omega$m. The background resistivity below the resistor is 2.5 $\Omega$m. To mimic a vertical transverse isotropic (VTI) Earth, all layers below the seabed are assigned with an anisotropy ratio (defined as $\lambda=\rho_v/\rho_h$ in this paper) $\lambda=1.5$. The source is placed in the middle of the model, 40 m above the receivers sitting on the seabed.
The resistivity in the water column and the formation above 1200 m was discretized with
constant grid spacing: $\Delta x=\Delta y = 150$ m, $ \Delta z =40$ m.  From 1200 m down to the bottom of the resistivity model, the grid has been stretched with different
growing factors.

\begin{figure}
  \centering
  \includegraphics[width=\linewidth]{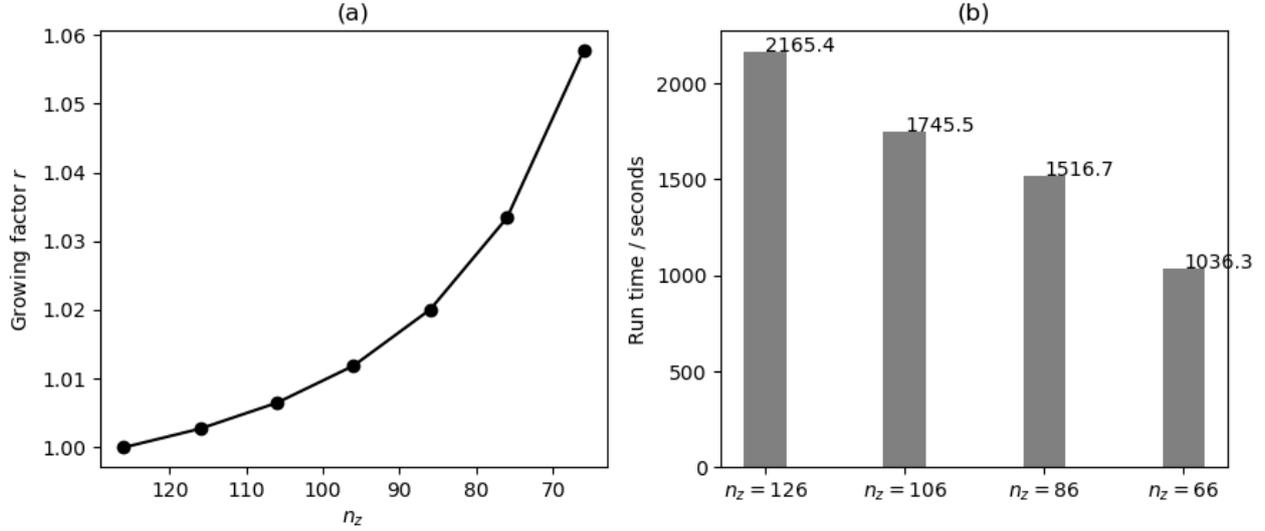}
  \caption{(a) The exponentially growing factor with the decreasing number of grid points $n_z$;
    (b) With the decreasing number of grid points $n_z$, the modelling time decreases dramatically.}\label{fig:figure7}
\end{figure}

Figure~\ref{fig:figure7}a displays the stretching factor $r$ of the nonuniform grid
for different number of grid points $n_z$. Note that the cell size grows very fast with
the factor $r$. For example, with $r=1.05$ after 40 cells we obtain $r^{40}>7$ times
the size of the first cell. Let us now increase $r$ (hence decrease $n_z$ correspondingly)
and analyze how it changes the total modelling time for the resistivity model of the
same physical length in $z$ direction. Figure~\ref{fig:figure7}b displays a significant
decrease of the modelling time associated with the reduction of the number of  grid points $n_z$.
We indeed see that modelling using $n_z=66$ takes almost half of the simulation time compared to modelling using $n_z=126$.

\begin{figure}
  \centering
  \includegraphics[width=\linewidth]{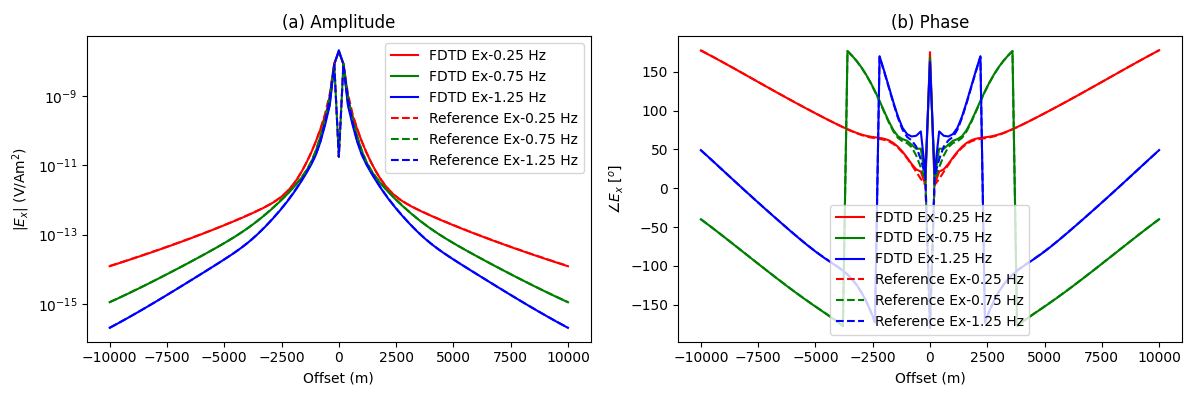}
  \caption{Comparison  of (a) amplitude and (b) phase between 4th order FDTD (solid line) and reference solution (dash line) for 3D CSEM simulation in the deep water scenario.}\label{fig:figure8}
\end{figure}

\begin{figure}
  \centering
  \includegraphics[width=\linewidth]{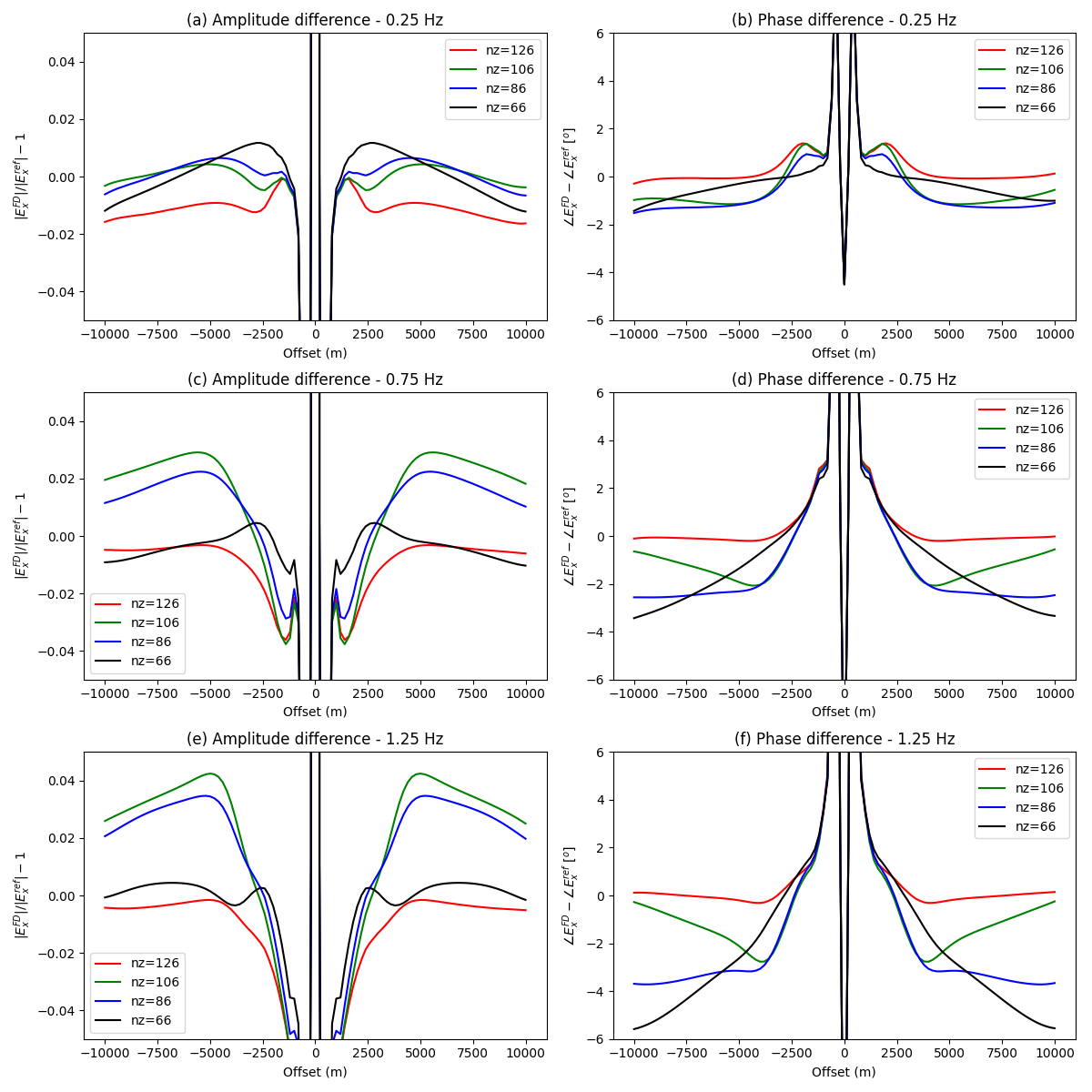}
  \caption{The amplitude and phase error of FDTD compared with 1D reference solution for 0.25 Hz (a,b), 0.75 Hz (c,d) and 1.25 Hz (e,f) in deep water.}\label{fig:figure9}
\end{figure}

Figure~\ref{fig:figure8} overlays the FDTD modelled EM fields with the reference solution. We see a very good agreement between the two in terms of both the amplitude and phase. 
To gain an idea of the modelling accuracy reduction when using less computing effort with nonuniform grid,
we plotted the amplitude and phase errors in Figure~\ref{fig:figure9} corresponding to three different frequencies.
These figures clearly show that both the amplitude error and the phase error increase with the increase of the frequencies. It is also interesting to note that stretching the grid does not necessarily increase the amplitude error, but does increase the phase discrepancy. This highlights the importance of combination in examing both amplitude and phase.

\subsection{3D modelling with seafloor bathymetry}

\begin{figure}
  \centering
  \includegraphics[width=0.8\linewidth]{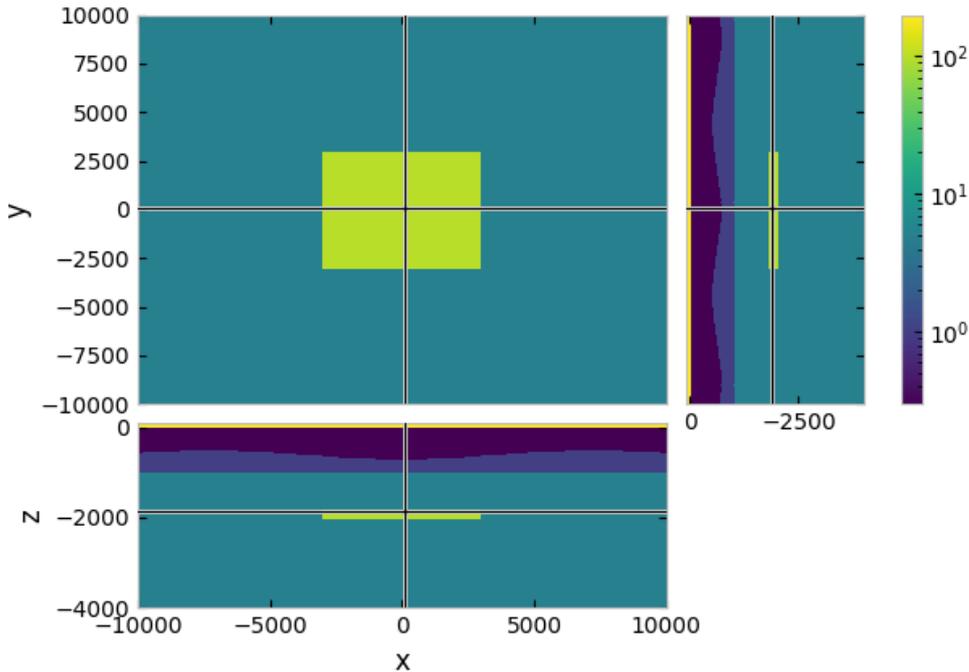}
  \caption{The 3D resistivity model with seafloor bathymetry, followed by 2 formation layers of resistivity values 1 $\Omega$m and 2 $\Omega$m. The last sediment layer includes a strong resistor of 100 $\Omega$m mimicking a hydrocarbon bearing formation located within the depth range [1800, 2000] m and the offset range [-3000, 3000] m in both x and y directions. The source was placed at the center of the model while 101 receivers with 200 m spacing are deployed along x direction.}\label{fig:figure10}
\end{figure}

The above 1D example highlights the importance of nonuniform grid in combination with
high-order schemes to achieve efficient numerical modelling with sufficient accuracy.
Here we consider a more realistic 3D resistivity model with bathymetry variations in horizontal directions, as shown in Figure~\ref{fig:figure10}a.

The model has seawater of 0.3 $\Omega$m, followed by 2 formation layers of resistivity - 1 $\Omega$m and 2 $\Omega$m. A resistor of 50 $\Omega$m was buried in the last formation layer to mimic a hydrocarbon bearing formation located between 1800 m and 2000 m in depth, with the offset expanding from -3000 m to 3000 m in both x and y directions.

In order to validate the modelling accuracy of the proposed method, we simulated a
reference solution by finite integration method in frequency domain using \verb|emg3d| \citep{werthmuller2019emg3d}
The finite-integration modelling extends the model tens of kilometers in each direction, to
mimic that the EM fields propagate to very far distance while avoiding possible edge/boundary
effects.  In our finite-difference modelling, the PML boundary condition attenuates the artificial
reflections in the computational domain within ten grid nodes to achieve the same behavior.
A dipole source was placed at 650 m depth in the middle of the model, while 101 receivers are sitting on the curved seabed.

\begin{figure}
  \centering
  \includegraphics[width=\textwidth]{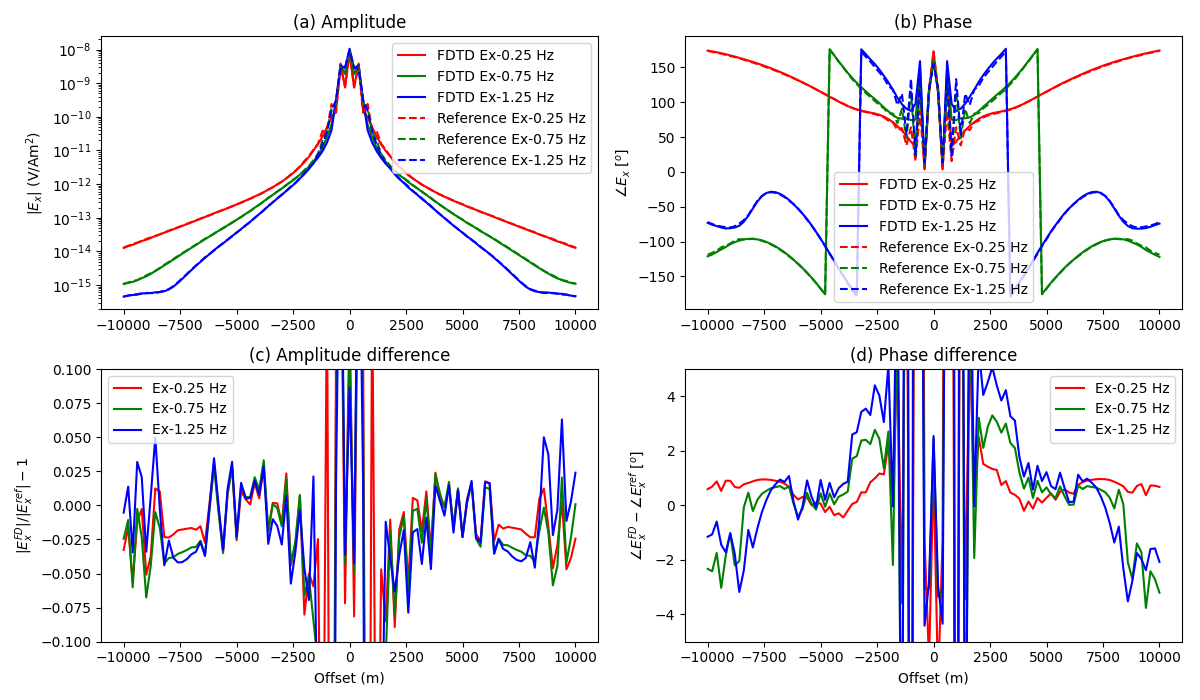}
  \caption{Comparison of the modelling results between the proposed method and the reference solution.}\label{fig:figure11}
\end{figure}

Figure~\ref{fig:figure11} displays the comparison of the 3D CSEM modelling between
our method and the result from \verb|emg3d|. From the
Figures~\ref{fig:figure11}a and \ref{fig:figure11}b, both the amplitude and
the phase from our method match very well with the reference solution. The maximum amplitude discrepancy for all frequencies in Figure~\ref{fig:figure11}c is bounded within 5\% at most of the relevant offset.
The maximum mismatch in phase is less than 1 degree for 0.25 Hz and 2 degrees for 0.75 Hz and 1.25 Hz for the EM field above the noise level (1e-15 V/m$^2$), as can be seen in
Figure~\ref{fig:figure11}d. These demonstrate the good accuracy achieved by our method.

\section{Discussion}

A natural idea to achieve higher modelling accuracy is to use denser sampling. However,
increasing the number of grid points in each dimension will lead to exponentially
growing computational overhead, i.e., double the sampling in x, y and z coordinates
results in eight times more nodes in the simulation. Over the same domain the node
separations are reduced with a factor of two. By that, the stability criterion dictates
a reduction in the time step by a factor of two. The net result is an increase of
computational cost with a factor of sixteen.
This is significantly more costly than considering
high-order FDTD.  For a fixed error requirement, high-order FDTD has been demonstrated
to be much more efficient than a lower-order scheme with dense sampling \citep{yefet2001staggered}.
Assuming a linear scaling dependency between computing time and the number of grid points,
doubling the length of the finite-difference operator in x, y and z directions
will simply double the computational cost.

Stretching along $z$-direction seems always to be beneficial in terms of computational
efficiency. It reduces the computational cost in two ways: first, it leads to increased grid
spacing, hence less grid points and larger node spacing to discretize the resistivity model for
simulation on the same physical size. Meanwhile, larger  grid spacing permits to use a larger
temporal step in terms of stability condition. This decreases the number of time steps needed to reach the steady state of the
frequency-domain EM fields.

It is natural to stretch the nonuniform grid in all three spatial directions, with
the motivation to decrease the computational cost further. In case there is no a priori
knowledge about the subsurface, a possible practice is to
start stretching from a given offset from the source location in the horizontal directions.
This approach is not followed up here due to the fact that grid stretching complicates
the calculation of the air-wave for time-domain codes. Grid stretching is applicable
also for frequency-domain finite-difference codes. The airwave implementation is very different
for frequency-domain codes, where the air layer is part of the simulation domain.
Frequency-domain codes lend themselves easily to both vertical and horizontal grid stretching.

\begin{figure}
  \centering
  \includegraphics[width=0.9\textwidth]{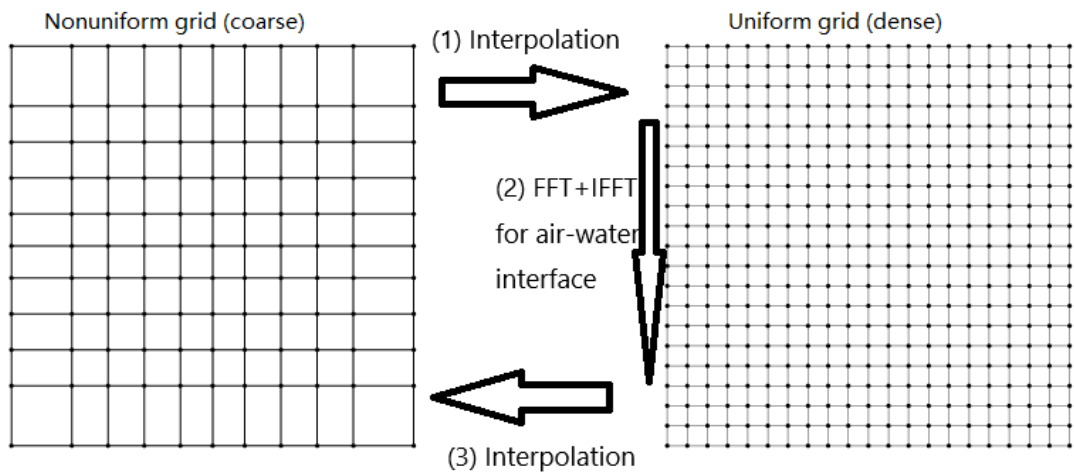}
  \caption{Horizontal grid stretching requires interpolating between the coarse
    nonuniform grid and a dense uniform grid, due to the equidistance
    requirement of FFT in airwave manipulation.}\label{fig:figure12}
\end{figure}

\begin{figure}
  \centering
  \includegraphics[width=\linewidth]{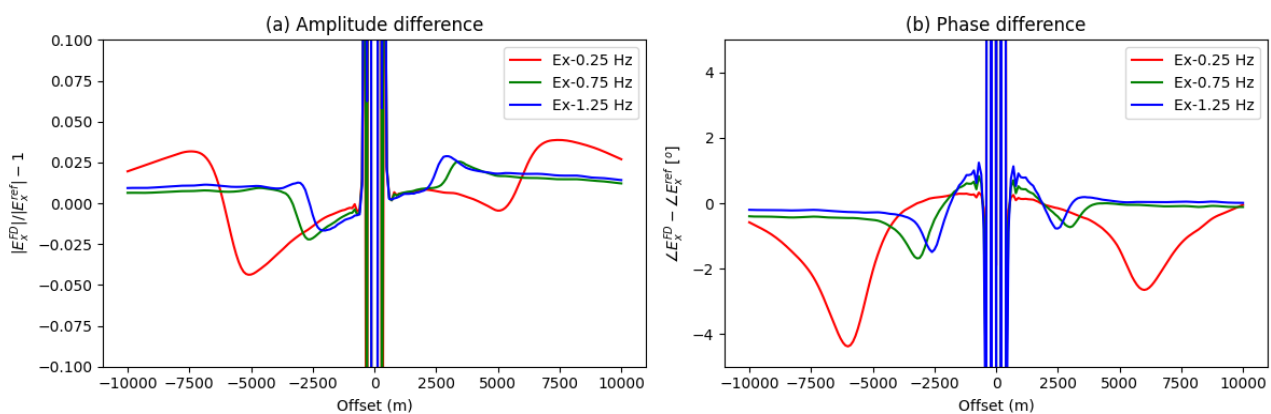}
  \caption{The amplitude and phase error after horizontal stretching of the 1D
    model with a stretching factor 1.05 along both $x$- and $y$- directions.}\label{fig:figure13}
\end{figure}

More research must be invested with respect to increasing the simulation efficiency
of time-domain finite-difference schemes focusing on airwave implementation combined
with horizontally nonuniform grids. The root of the problem is that the
most common airwave implementations require fast Fourier transform (FFT) on
regular grids. The fields are transformed to the wavenumber domain, propagated into the
air layer in this domain and then transformed back to the space domain.
A straight forward approach is to use interpolation  between the uniform grid
and the nonuniform one, as illustrated in Figure~\ref{fig:figure12}.
The grid stretching in $x$- and $y$- directions certainly complicates the implementation
while introducing additional computational cost, which is opposite to what we want
to achieve. We have tested horizontal stretching for the previous 1D model
using the stretching factor 1.05 in both $x$- and $y$- directions. The resulting
amplitude and phase error in Figure~\ref{fig:figure13} shows that the numerical
accuracy is highly degraded (maximum amplitude error for 0.25 Hz is around 4\%)
compared with the result for a horizontally uniform grid
(maximum amplitude error is less than 1.5\% for all calculated frequencies).
The running time is in fact longer than for the horizontally uniform grid which has
a higher number of nodes. Due to the nonuniform grid staggering, several nodes on the
uniform grid may reside in the same interval between two neighboring nodes on the
nonuniform grid. The error panel displays an unsymmetrical pattern
in Figure~\ref{fig:figure13}. Interpolating the fields from the uniform grid
back to nonuniform grid can thus produce less accurate solutions. No stretching
is therefore recommended along $x$- and $y$- for efficiency and accuracy considerations
until a more accurate solution to this problem is developed.

 It is noteworthy that fictitious wave domain is simply a mathematical tool to compute correct frequency domain EM response efficiently by time stepping. The fictitious time is different from real time. To compute the true EM time series correctly, one needs to simulate a large number frequencies of the CSEM fields and then perform inverse Fourier transform, as has been done in \citet{Mittet_2010_HFD} and \citet{rochlitz2021evaluation}. 

\section{Conclusion}

We have presented a 3D CSEM modelling method using high-order FDTD on a nonuniform grid.
The key problem addressed in this work is the low
accuracy and inconsistency issue in standard 2nd order staggered FDTD scheme on nonuniform
grid. The strategy we propose is to adapt the interpolation weights depending on the nodal
distance. These finite difference coefficients can be computed by inverting a Vandermonde matrix in an accurate and efficient manner. This makes our approach different from the commonly used EM modelling approaches. A new yet more generic stability condition has been established in order to achieve stable FDTD modelling.
In designing the nonuniform grid based on geometrical progression, we develop a fixed point iteration to compute the
optimal growing factor which allows good match of the modelling domain in case of grid
stretching. The numerical examples demonstrate that there
is a significant improvement in accuracy by using a high-order FDTD scheme, while
combining it with a nonuniform grid reduces the computational cost without a significant
sacrifice of accuracy. We conclude that high-order finite differences on nonuniform
grid is a viable tool for full scale 3D CSEM modelling applications.
Since the key idea is to use high order finite-difference coefficients adaptive to the node spacing, the method is expected to be applicable also for finite-difference frequency-domain schemes.

\section*{Acknowledgments}

Pengliang Yang was supported by Chinese Fundamental Research Funds for the Central
Universities (AUGA5710010121) and National Natural Science Fundation of China (42274156). Pengliang Yang thanks Dieter Werthermullter for
the assistance to produce the reference solution  using \verb|empymod| and \verb|emg3d| to validate the accuracy of the proposed method. The source code of this work can be found in the github repository:
\url{https://github.com/yangpl/libEMM}.

\appendix

\section*{Computing uniform staggered-grid finite difference coefficients via Vandermonde matrix inversion}\label{appendix:fdcoeff}

The method to invert Vandermonde matrices gives us a generic approach to compute finite-difference
coefficients with arbitrary grid spacing. A special case is the regular grid spacing. In what follows,
we show how the standard staggered grid finite difference coefficients can be accurately computed also
within the same framework.

The Taylor series expansion of a function $f(x)$ can be written as
\begin{equation}\label{eq:Taylor}
  \begin{cases}
    f(x+h)=f(x)+\frac{\partial f(x)}{\partial x}h+\frac{1}{2!}\frac{\partial^2 f(x)}{\partial x^2}h^2+\frac{1}{3!}\frac{\partial^3 f(x)}{\partial x^3}h^3+\ldots\\
    f(x-h)=f(x)-\frac{\partial f(x)}{\partial x}h+\frac{1}{2!}\frac{\partial^2 f(x)}{\partial x^2}h^2-\frac{1}{3!}\frac{\partial^3 f(x)}{\partial x^3}h^3+\ldots
  \end{cases}.
\end{equation}
It leads to
\begin{equation}
  \begin{cases}
    \frac{f(x+h)+f(x-h)}{2}
    &=f(x)+\frac{1}{2!}\frac{\partial^2 f(x)}{\partial x^2}h^2+\frac{1}{4!}\frac{\partial^4 f(x)}{\partial x^4}h^4+\ldots\\
    \frac{f(x+h)-f(x-h)}{2}
    &=\frac{\partial f(x)}{\partial x}h+\frac{1}{3!}\frac{\partial^3 f(x)}{\partial x^3}h^3+\frac{1}{5!}\frac{\partial^5 f(x)}{\partial x^5}h^5+\ldots
  \end{cases}.
\end{equation}
Let $h=\Delta x/2$. This implies the 2nd order accuracy of centered finite difference scheme using only two staggered nodes
\begin{equation}\label{eq:approx}
  \begin{cases}
    \frac{\partial f(x)}{\partial x}=\frac{f(x+\Delta x/2)-f(x-\Delta x/2)}{\Delta x}+O(\Delta x^2)\\
    f(x)=\frac{f(x+\Delta x/2)+f(x-\Delta x/2)}{2}+O(\Delta x^2)
  \end{cases}.
\end{equation}

To approximate the 1st order derivatives as accurate as possible, we express it using more
consecutive nodes. Due to regular grid staggering, the coefficients lying on symmetric
positions should have the same coefficients. This means the first-order derivative reads in
the following form
\begin{equation}
  \begin{split}
    \frac{\partial f}{\partial x}=&b_1\frac{f(x+\Delta x/2)-f(x-\Delta x/2)}{\Delta x}+\\
    &b_2\frac{f(x+3\Delta x/2)-f(x-3\Delta x/2)}{3\Delta x}+\\
    &b_3\frac{f(x+5\Delta x/2)-f(x-5\Delta x/2)}{5\Delta x}+\cdots.
  \end{split}
\end{equation}
Substituting the $f(x+h)$ and $f(x-h)$ with
equation \ref{eq:Taylor} for $h=\Delta x/2,3\Delta x/2, \ldots$ results in
\begin{equation}
  \begin{split}
    \frac{\partial f}{\partial x}=&b_1\cdot 2\left(\frac{\Delta x}{2} \frac{\partial f}{\partial x}+\frac{1}{3!}(\frac{\Delta x}{2})^3\frac{\partial^3 f}{\partial x^3}+\cdots\right)/{\Delta x}\\
    &+b_2\cdot 2\left(\frac{3\Delta x}{2} \frac{\partial f}{\partial x}+\frac{1}{3!}(\frac{3\Delta x}{2})^3\frac{\partial^3 f}{\partial x^3}+\cdots\right)/{3\Delta x}\\
    &+b_3\cdot 2\left(\frac{5\Delta x}{2} \frac{\partial f}{\partial x}+\frac{1}{3!}(\frac{5\Delta x}{2})^3\frac{\partial^3 f}{\partial x^3}+\cdots\right)/{5\Delta x}+\ldots\\
    =&(b_1+b_2+b_3+b_4+\cdots)\frac{\partial f}{\partial x}\\
    &+\frac{\Delta x^2}{3!\cdot 2^2}(b_1+3^2b_2+5^2b_3+7^2b_4+\cdots)\frac{\partial^3 f}{\partial x^3}\\
    &+\frac{\Delta x^4}{5!\cdot 2^4}(b_1+3^4b_2+5^4b_3+7^4b_4+\cdots)\frac{\partial^5 f}{\partial x^5}+\cdots.
  \end{split}
\end{equation}
Thus, taking first $L$ terms (corresponding to using $2L$ nodes) requires
\begin{equation}
  \begin{cases}
    b_1+b_2+b_3+\cdots+b_L&=1\\
    b_1+3^2b_2+5^2b_3+\cdots+(2L-1)^2b_L&=0\\
    b_1+3^4b_2+5^4b_3+\cdots+(2L-1)^4b_L&=0\\
    \cdots	&\\
    b_1+3^{2L-2}b_2+5^{2L-2}b_3+\cdots+(2L-1)^{2L-2}b_L&=0\\
  \end{cases},
\end{equation}
which again builds up a Vandermonde-like system
\begin{gather}\label{eq:Vab}
  \underbrace{
    \begin{bmatrix}
      1 & 1 & \ldots & 1\\
      x_1 & x_2 & \ldots & x_L\\
      \vdots &  & \ddots & \vdots\\
      x_1^{L-1} & x_2^{L-1} & \ldots & x_L^{L-1}\\
    \end{bmatrix}
  }_{\textbf{V}}
  \underbrace{
    \begin{bmatrix}
      b_1\\
      b_2\\
      \vdots\\
      b_L\\
    \end{bmatrix}
  }_{\textbf{b}}=
  \underbrace{
    \begin{bmatrix}
      1\\
      0\\
      \vdots\\
      0\\
    \end{bmatrix}
  }_{\textbf{z}},
\end{gather}
in which $x_i=(2i-1)^2$, $i=1,\cdots,L$. The finite difference weights can then be easily
computed using \citet[Algorithm 4.6.2]{Golub_1996_MATCOMP}, see an Octave/Matlab script
for computing them in \citet[section 2.5]{Yang_2014_NTW}. These numerically computed
weights may be cross-validated with the generic method by \citet[Table 2]{Fornberg_1988_GFD}.
With regular grid spacing $\Delta x$,  the  weights $b_i$ are dimensionless and can be connected
to the coefficients $\alpha_i$ in equation \ref{eq:gandalf} involving a scaling factor $\Delta x$.

\bibliographystyle{seg}

\begin{thebibliography}{}
\itemsep0pt

\bibitem[Bj{\"o}rck and Pereyra, 1970]{bjorck1970solution}
Bj{\"o}rck, A., and V. Pereyra,  1970, {Solution of Vandermonde systems of
  equations}: Mathematics of computation, {\bfseries 24}, 893--903.

\bibitem[da~Silva et~al., 2012]{da2012finite}
da~Silva, N.~V., J.~V. Morgan, L. MacGregor, and M. Warner,  2012, {A finite
  element multifrontal method for {3D} CSEM modeling in the frequency domain}:
  Geophysics, {\bfseries 77}, E101--E115.

\bibitem[de~Hoop, 1996]{de1996general}
de~Hoop, A.~T.,  1996, A general correspondence principle for time-domain
  electromagnetic wave and diffusion fields: Geophysical Journal International,
  {\bfseries 127}, 757--761.

\bibitem[Demmel and Koev, 2005]{demmel2005accurate}
Demmel, J., and P. Koev,  2005, The accurate and efficient solution of a
  totally positive generalized vandermonde linear system: SIAM Journal on
  Matrix Analysis and Applications, {\bfseries 27}, 142--152.

\bibitem[Fornberg, 1988]{Fornberg_1988_GFD}
Fornberg, B.,  1988, Generation of finite difference formulas on arbitrarily
  spaced grids: Mathematics of Computation, {\bfseries 51}, 699--706.

\bibitem[Golub, 1996]{Golub_1996_MATCOMP}
Golub, G.~H.,  1996, Matrix computation, third edition: Johns Hopkins Studies
  in Mathematical Sciences.

\bibitem[Key, 2016]{key2016mare2dem}
Key, K.,  2016, {MARE2DEM: a 2-D inversion code for controlled-source
  electromagnetic and magnetotelluric data}: Geophysical Journal International,
  {\bfseries 207}, 571--588.

\bibitem[Komatitsch and Martin, 2007]{Komatitsch_2007_GEO}
Komatitsch, D., and R. Martin,  2007, {An unsplit convolutional perfectly
  matched layer improved at grazing incidence for the seismic wave equation}:
  Geophysics, {\bfseries 72}, SM155--SM167.

\bibitem[Lee et~al., 1989]{lee1989new}
Lee, K.~H., G. Liu, and H. Morrison,  1989, A new approach to modeling the
  electromagnetic response of conductive media: Geophysics, {\bfseries 54},
  1180--1192.

\bibitem[Li and Key, 2007]{li20072d}
Li, Y., and K. Key,  2007, {2D} marine controlled-source electromagnetic
  modeling: Part 1—an adaptive finite-element algorithm: Geophysics,
  {\bfseries 72}, WA51--WA62.

\bibitem[Maa{\o}, 2007]{Maao_2007_FFT}
Maa{\o}, F.,  2007, Fast finite-difference time-domain modeling for marine
  subsurface electromagnetic problems: Geophysics, {\bfseries 72}, A19--A23.

\bibitem[Mittet, 2010]{Mittet_2010_HFD}
Mittet, R.,  2010, High-order finite-difference simulations of marine {CSEM}
  surveys using a correspondence principle for wave and diffusion fields:
  Geophysics, {\bfseries 75}, F33--F50.

\bibitem[Mittet, 2015]{mittet2015seismic}
--------, 2015, Seismic wave propagation concepts applied to the interpretation
  of marine controlled-source electromagnetics: Geophysics, {\bfseries 80},
  E63--E81.

\bibitem[Mittet, 2017]{mittet2017internal}
--------, 2017, On the internal interfaces in finite-difference schemes:
  Geophysics, {\bfseries 82}, T159--T182.

\bibitem[Monk and S{\"u}li, 1994]{monk1994convergence}
Monk, P., and E. S{\"u}li,  1994, A convergence analysis of {Y}ee’s scheme on
  nonuniform grids: SIAM Journal on Numerical Analysis, {\bfseries 31},
  393--412.

\bibitem[Mulder, 2006]{mulder2006multigrid}
Mulder, W.,  2006, A multigrid solver for 3{D} electromagnetic diffusion:
  Geophysical prospecting, {\bfseries 54}, 633--649.

\bibitem[Newman and Alumbaugh, 1995]{newman1995frequency}
Newman, G.~A., and D.~L. Alumbaugh,  1995, Frequency-domain modelling of
  airborne electromagnetic responses using staggered finite differences:
  Geophysical Prospecting, {\bfseries 43}, 1021--1042.

\bibitem[Oristaglio and Hohmann, 1984]{oristaglio1984diffusion}
Oristaglio, M.~L., and G.~W. Hohmann,  1984, Diffusion of electromagnetic
  fields into a two-dimensional earth: A finite-difference approach:
  Geophysics, {\bfseries 49}, 870--894.

\bibitem[Puzyrev et~al., 2013]{puzyrev2013parallel}
Puzyrev, V., J. Koldan, J. de~la Puente, G. Houzeaux, M. V{\'a}zquez, and J.~M.
  Cela,  2013, A parallel finite-element method for three-dimensional
  controlled-source electromagnetic forward modelling: Geophysical Journal
  International, {\bfseries 193}, 678--693.

\bibitem[Rochlitz et~al., 2021]{rochlitz2021evaluation}
Rochlitz, R., M. Seidel, and R.-U. Börner,  2021, {Evaluation of three
  approaches for simulating 3-D time-domain electromagnetic data}: Geophysical
  Journal International, {\bfseries 227}, 1980--1995.

\bibitem[Rochlitz et~al., 2019]{rochlitz2019custem}
Rochlitz, R., N. Skibbe, and T. G{\"u}nther,  2019, custem: Customizable
  finite-element simulation of complex controlled-source electromagnetic data:
  Geophysics, {\bfseries 84}, F17--F33.

\bibitem[Roden and Gedney, 2000]{roden2000convolution}
Roden, J.~A., and S.~D. Gedney,  2000, {Convolution PML (CPML): An efficient
  FDTD implementation of the CFS--PML for arbitrary media}: Microwave and
  optical technology letters, {\bfseries 27}, 334--339.

\bibitem[Smith, 1996a]{smith1996conservative1}
Smith, J.~T.,  1996a, Conservative modeling of 3-{D} electromagnetic fields,
  {P}art i: Properties and error analysis: Geophysics, {\bfseries 61},
  1308--1318.

\bibitem[Smith, 1996b]{smith1996conservative2}
--------, 1996b, {Conservative modeling of 3-D electromagnetic fields, Part II:
  Biconjugate gradient solution and an accelerator}: Geophysics, {\bfseries
  61}, 1319--1324.

\bibitem[Streich, 2009]{streich20093d}
Streich, R.,  2009, 3{D} finite-difference frequency-domain modeling of
  controlled-source electromagnetic data: {D}irect solution and optimization
  for high accuracy: Geophysics, {\bfseries 74}, F95--F105.

\bibitem[Taflove and Hagness, 2005]{Taflove_2005_CEF}
Taflove, A., and S.~C. Hagness,  2005, Computational electrodynamics: The
  finite-difference time-domain method, 3rd ed.: Artech House.

\bibitem[Wang and Hohmann, 1993]{wang1993finite}
Wang, T., and G.~W. Hohmann,  1993, A finite-difference, time-domain solution
  for three-dimensional electromagnetic modeling: Geophysics, {\bfseries 58},
  797--809.

\bibitem[Werthmüller, 2017]{werthmuller2017open}
Werthmüller, D.,  2017, {An open-source full 3D electromagnetic modeler for 1D
  VTI media in Python: empymod}: Geophysics, {\bfseries 82}, WB9--WB19.

\bibitem[Werthmüller et~al., 2019]{werthmuller2019emg3d}
Werthmüller, D., W. Mulder, and E. Slob,  2019, emg3d: A multigrid solver for
  3d electromagnetic diffusion: Journal of Open Source Software, {\bfseries 4},
  1463.

\bibitem[Yang, 2014]{Yang_2014_NTW}
Yang, P.,  2014, A numerical tour of wave propagation: Technical report, Xi'an
  Jiaotong University.

\bibitem[Yee, 1966]{Yee_1966_NSI}
Yee, K.~S.,  1966, Numerical solution of initial boundary value problems
  involving {M}axwell's equations in isotropic media: {IEEE} Transactions on
  Antennas and Propagation, {\bfseries 14}, 302--307.

\bibitem[Yefet and Petropoulos, 2001]{yefet2001staggered}
Yefet, A., and P.~G. Petropoulos,  2001, A staggered fourth-order accurate
  explicit finite difference scheme for the time-domain maxwell's equations:
  Journal of Computational Physics, {\bfseries 168}, 286--315.

\end{thebibliography}

\newcommand{\SortNoop}[1]{}

\end{document}